\documentclass[review,hidelinks,onefignum,onetabnum]{siamart250211}

\usepackage{lipsum}
\usepackage{amsfonts}
\usepackage{graphicx}
\usepackage{epstopdf}
\usepackage{algpseudocode}
\usepackage{setspace}
\usepackage{subcaption}
\let\Algorithm\algorithm
\renewcommand\algorithm[1][]{\Algorithm[#1]\setstretch{1.15}}

\ifpdf
  \DeclareGraphicsExtensions{.eps,.pdf,.png,.jpg}
\else
  \DeclareGraphicsExtensions{.eps}
\fi

\newsiamremark{remark}{Remark}
\newsiamremark{hypothesis}{Hypothesis}
\crefname{hypothesis}{Hypothesis}{Hypotheses}
\newsiamthm{claim}{Claim}
\newsiamremark{fact}{Fact}
\crefname{fact}{Fact}{Facts}

\headers{On efficient block Krylov solvers for $\mathcal H^2$-matrices}{S. Christophersen}

\title{On efficient block Krylov solvers for $\mathcal H^2$-matrices}

\author{Sven Christophersen \thanks{Mathematisches Seminar, Christian-Albrechts-Universität Kiel,
        Heinrich-Hecht-Platz 6, Kiel, Germany, \email{christophersen@math.uni-kiel.de}}}

\usepackage{amsopn}

\ifpdf
\hypersetup{
  pdftitle={On efficient block Krylov-solvers for $\mathcal H^2$-matrices},
  pdfauthor={S. Christophersen}
}
\fi

\newcommand{\Idx}{\mathcal I}
\newcommand{\Jdx}{\mathcal J}
\newcommand{\Tree}{\mathcal T}
\newcommand{\Leafs}{\mathcal L}
\newcommand{\TI}{\Tree_\Idx}
\newcommand{\TJ}{\Tree_\Jdx}
\newcommand{\TIJ}{\Tree_{\Idx \times \Jdx}}
\newcommand{\LI}{\Leafs_\Idx}

\newcommand{\LIJ}{\Leafs_{\Idx \times \Jdx}}
\newcommand{\Clist}{\mathcal C}

\newcommand{\sons}{\operatorname{sons}}
\newcommand{\Root}{\operatorname{root}}
\newcommand{\diam}{\operatorname{diam}}
\newcommand{\dist}{\operatorname{dist}}
\newcommand{\supp}{\operatorname{supp}}

\algblock{ParFor}{EndParFor}
\algnewcommand\algorithmicparfor{\textbf{parfor}}
\algnewcommand\algorithmicpardo{\textbf{do}}
\algnewcommand\algorithmicendparfor{\textbf{end\ parfor}}
\algrenewtext{ParFor}[1]{\algorithmicparfor\ #1\ \algorithmicpardo}
\algrenewtext{EndParFor}{\algorithmicendparfor}

\begin{document}
\nolinenumbers

\maketitle

% REQUIRED
\begin{abstract}
Hierarchical matrices provide a highly memory-efficient way of storing dense linear operators
arising, for example, from boundary element methods, particularly when stored in the 
$\mathcal H^2$ format. In such data-sparse representations, iterative solvers are preferred 
over direct ones due to the cost-efficient matrix–vector multiplications they enable.
Solving multiple systems of linear equations with the same hierarchical matrix naturally 
leads to block methods, which in turn make heavy use of BLAS level-3 functions such as \texttt{GEMM}. \\
We present an efficient implementation of $\mathcal H^2$-matrix–vector and 
$\mathcal H^2$-matrix–matrix multiplication that fully exploits the potential of modern 
hardware in terms of memory and cache utilization. The latter is employed to accelerate 
block Krylov subspace methods, which we present later as the main results of this paper.
\end{abstract}

% REQUIRED

% REQUIRED

\section{Introduction}

In almost every numerical application, the task of solving a linear system of equations arises naturally. \\ 
Let $n \in \mathbb N$ and let $\mathbf{A} \in \mathbb K^{n \times n}$ 
be an invertible matrix. Furthermore, let $b \in \mathbb K^n$ denote the 
right-hand side and $x \in \mathbb K^n$ the solution to the system
\[
\mathbf{A} x = b .
\]
Depending on the structure of the matrix $\mathbf{A}$, different techniques 
are employed to solve such systems. For dense matrices, one usually prefers 
matrix factorizations such as LU or Cholesky decompositions, followed by 
forward and backward substitution. This approach typically has cubic 
complexity with respect to $n$ for the factorization and $\mathcal O(n^2)$
complexity for the subsequent solve operations and is therefore not applicable
to large systems.
On the other hand, when the structure of the matrix is favorable and the
memory footprint is low, iterative solvers such as Krylov subspace methods
are often applied to compute the solution.

When dealing with integral equations, a particularly successful approach to 
handle the resulting dense matrices is the use of \emph{hierarchical matrices} 
in various forms, such as $\mathcal H$-matrices or more advanced variants like 
$\mathcal H^2$- or $\mathcal{DH}^2$-matrices. 
The central idea is to compress certain submatrices of the matrix by 
low-rank approximations of rank at most $k \in \mathbb N$. 
This technique significantly reduces both the storage requirements and the 
computational complexity of the matrix–vector product, which is 
$\mathcal O(n k \log n)$ or even $\mathcal O(n k)$.

Since the convergence rate of Krylov methods such as \emph{CG} ~\cite{hestenes1952methods}
or \emph{GMRES} ~\cite{saad1986gmres} strongly depends on the spectrum of $\mathbf{A}$, an additional preconditioner 
is often constructed, for example as an approximate LU or Cholesky decomposition
in the $\mathcal H$-matrix format, leading to significantly reduced iteration
counts ~\cite{bebendorf2005hierarchical, grasedyck2003construction, aliaga2017task, carratala2019exploiting}.
This approach is well established and performs reliably in practical applications.

In the course of this article, we consider a slightly different setting: 
we are not only interested in solving a single linear system, but in computing the 
solutions of $m \in \mathbb N$ systems, i.e., we aim to solve
\begin{equation}
\mathbf A x^{(i)} = b^{(i)}, \qquad \text{for } i = 1, \ldots , m .
\end{equation}
Accordingly, we introduce the matrices 
\[
\mathbf{x} := 
\begin{pmatrix} x^{(1)} & \cdots & x^{(m)} \end{pmatrix} \in \mathbb{K}^{n \times m},
\qquad
\mathbf{b} := 
\begin{pmatrix} b^{(1)} & \cdots & b^{(m)} \end{pmatrix} \in \mathbb{K}^{n \times m},
\]
and reformulate our problem in terms of matrix multiplication as
\begin{equation}\label{eqn:matrix_system}
\mathbf{A} \mathbf{x} = \mathbf{b} .
\end{equation}
We will demonstrate how to implement equation~\eqref{eqn:matrix_system} in the 
context of Krylov methods, making use of BLAS level-3 routines such as 
\texttt{GEMM} instead of the standard matrix–vector products, such as
\texttt{GEMV}, which belong to the domain of BLAS level-2.

To achieve a high utilization of the available hardware, 
we propose a novel approach to parallelizing the matrix–vector product of an 
$\mathcal H^2$-matrix with a dense vector, as well as the corresponding product 
with a dense matrix, which naturally arises in the context of block Krylov 
subspace methods discussed later.

Our focus lies exclusively on exploiting BLAS level-3 routines 
to optimize hardware efficiency in the context of $\mathcal H^2$-matrices. 
We do not seek to enrich the Krylov subspace itself, as has been done in previous 
works ~\cite{gutknecht2005block, gutknecht2007block, gutknecht2009block, rashedi2016short, o1980block, morgan2005restarted, simoncini1995iterative, simoncini1996hybrid}.

The key novelty of this paper is the combination of hierarchical matrix techniques 
with blocked Krylov subspace methods, yielding highly efficient solvers for 
multiple linear systems that share the same matrix but have independent right-hand sides.

\section{Hierarchical matrices}
\label{sec:hmatrices}

Before turning our attention to solving the linear systems described above, 
we briefly recall the basic concept of hierarchical matrices
~\cite{Hackbusch2002, hackbusch2015hierarchical, bebendorf2008hierarchical, borm2010efficient}.

These matrices are frequently applied in the context of boundary element methods, 
which we use here as a motivating example. 
Let $\Omega \subseteq \mathbb R^3$ be a domain, and for general index sets 
$\Idx$ and $\Jdx$ consider families of basis functions $(\varphi_i)_{i \in \Idx}$ 
and $(\psi_j)_{j \in \Jdx}$ defined on $\Omega$. 
Furthermore, let $g : \Omega \times \Omega \to \mathbb K$ be a kernel function. 
We are interested in computing the matrix $G \in \mathbb K^{\Idx \times \Jdx}$ defined by
\begin{equation}
G_{ij} = \int_\Omega \int_\Omega \varphi_i(x)\, g(x,y)\, \psi_j(y)\, dy\, dx, 
\qquad \text{for all } i \in \Idx, j \in \Jdx .
\label{eqn:integral_eq}
\end{equation}
The resulting matrix is densely populated and therefore requires 
$\mathcal O(n^2)$ units of storage if $|\Idx| = |\Jdx| = n$. 
Hence, we seek a data-sparse approximation of $G$ in order to handle large-scale 
problems efficiently. \\

The idea is to identify suitable subsets $\hat t \subseteq \Idx$ and 
$\hat s \subseteq \Jdx$ such that the subdomains
\[
\Omega_t := \bigcup_{i \in \hat t} \supp(\varphi_i) \subseteq \Omega,
\qquad 
\Omega_s := \bigcup_{j \in \hat s} \supp(\psi_j) \subseteq \Omega
\]
satisfy the condition that $g$ is smooth on at least $\Omega_t \times \Omega_s$ 
and can therefore be approximated by a low-rank matrix.

This means, the task is to identify such subsets of the index sets $\Idx$ and $\Jdx$, 
which can be achieved by means of so-called \emph{cluster trees}.

\begin{definition}[Cluster tree]
Let $\Idx$ be a finite index set and let $\Tree$ be a labeled tree.
For a node $t \in \Tree$, we denote its \emph{label} by $\hat t$.
The tree $\Tree$ is called a \emph{cluster tree} for the index set $\Idx$ if
\begin{itemize}
\item the label of the root $r = \Root(\Tree)$ satisfies $\hat r = \Idx$, 
\item for every $t \in \Tree$ with $\sons(t) \neq \emptyset$, we have
      $\hat t = \bigcup_{t' \in \sons(t)} \hat{t'}$, and
\item for every $t \in \Tree$ and for all $t_1, t_2 \in \sons(t)$ with $t_1 \neq t_2$, 
      the labels are disjoint, i.e.\ $\hat{t_1} \cap \hat{t_2} = \emptyset$.
\end{itemize}
To avoid ambiguity, we often denote the tree by $\TI$.
Its nodes are referred to as \emph{clusters}, and the set of leaves is defined by
\[
\LI := \{\, t \in \TI : \sons(t) = \emptyset \,\}.
\]
\end{definition}

The construction of cluster trees is performed recursively by splitting a given 
index set into two or more non-overlapping parts. 
This process is repeated until each cluster contains a minimal number of elements, 
at which point the recursion terminates.

Next, we aim to construct a disjoint partition of the index set $\Idx \times \Jdx$, 
which is achieved by a so-called \emph{block tree}.

\begin{definition}[Block tree]\label{def:block_tree}
Let $\TI$ and $\TJ$ be cluster trees for the sets $\Idx$ and $\Jdx$, respectively. 
The tree $\TIJ$ is called a \emph{block tree} for $\TI$ and $\TJ$ if
\begin{itemize}
\item for every $b \in \TIJ$, there exist $t \in \TI$ and $s \in \TJ$ such that 
      $b = (t,s)$ and $\hat b = \hat t \times \hat s$,
\item for the root $b = \Root(\TIJ)$, we have $b = (\Root(\TI), \Root(\TJ))$,
\item for every $b = (t,s) \in \TIJ$ with $\sons(b) \neq \emptyset$, 
      the children of $b$ are given by
\[
\sons(b) =
\begin{cases}
  \{ t \} \times \sons(s) &: \text{if } \sons(t) = \emptyset , \\
  \sons(t) \times \{ s \} &: \text{if } \sons(s) = \emptyset , \\
  \sons(t) \times \sons(s) &: \text{otherwise} .
\end{cases}
\]
\end{itemize}
The nodes of $\TIJ$ are called \emph{blocks}, and, analogously to cluster trees, 
the set of leaves is denoted by $\LIJ$.
\end{definition}

Since the set $\LIJ$ defines a disjoint partition of the entire index set 
$\Idx \times \Jdx$, we can use the leaves of $\TIJ$ to efficiently handle 
sub-blocks of matrices $G \in \mathbb K^{\Idx \times \Jdx}$.

The final ingredient is a suitable \emph{admissibility condition}, for instance 
the standard condition
\begin{equation}\label{eqn:admiss}
\max \{ \diam(\Omega_t), \diam(\Omega_s) \} \leq \eta \, \dist(\Omega_t, \Omega_s),
\end{equation}
for some $\eta > 0$, which evaluates to either true or false for given subsets 
$\hat t \subseteq \Idx$ and $\hat s \subseteq \Jdx$. 
We then distinguish between admissible leaves 
\[
\LIJ^{+} := \{ b \in \LIJ : \eqref{eqn:admiss} \text{ holds} \}
\]
and inadmissible leaves 
\[
\LIJ^{-} := \LIJ \setminus \LIJ^{+}.
\] 
For the latter, the corresponding matrix entries are stored directly, i.e., 
$G_{\hat t \times \hat s}$ is computed and stored explicitly.

For admissible blocks $b = (t,s) \in \LIJ^{+}$, we can find matrices 
$V_t \in \mathbb K^{\hat t \times k}$ and $W_s \in \mathbb K^{\hat s \times k}$, 
as well as a small matrix $S_b \in \mathbb K^{k \times k}$, such that the block 
admits a rank-$k$ approximation
\begin{equation}\label{eqn:h2_admissible}
G_{\hat t \times \hat s} \approx V_t \, S_b \, W_s^* .
\end{equation}
The small matrices $S_b$ for $b \in \LIJ^{+}$ are denoted as \emph{coupling
matrices} and the families $(V_t)_{t \in \TI}$ and $(W_s)_{s \in \TJ}$ are
referred to as \emph{cluster bases}, and we will explain their construction in
the following.

\begin{definition}[Cluster basis]
Let $k \in \mathbb N$ be the local rank, and let $V = (V_t)_{t \in \TI}$ 
be a family of matrices with $V_t \in \mathbb K^{\hat t \times k}$ for all $t \in \TI$. 
If for all $t \in \TI \setminus \LI$ and all $t' \in \sons(t)$ there exist matrices 
$E_{t'} \in \mathbb K^{k \times k}$ such that
\[
V_t|_{\hat{t'} \times k} = V_{t'} \, E_{t'},
\]
then the family $V$ is called a \emph{(nested) cluster basis}. 
The matrices $(E_t)_{t \in \TI}$ are referred to as \emph{transfer matrices}, 
and the matrices $(V_t)_{t \in \LI}$ as \emph{leaf matrices}.
\end{definition}

Next, we introduce the final definition of what we mean by an \emph{$\mathcal H^2$-matrix}.

\begin{definition}[$\mathcal H^2$-matrix]
Let $G \in \mathbb K^{\Idx \times \Jdx}$ be a matrix, and let $\TIJ$ be an admissible 
block tree. Further, let $V = (V_t)_{t \in \TI}$ and $W = (W_s)_{s \in \TJ}$ 
be nested cluster bases. 
If the representation \eqref{eqn:h2_admissible} holds for all pairs 
$b = (t,s) \in \LIJ^+$, then $G$ is called an $\mathcal H^2$-matrix of rank $k \in \mathbb N$.
\end{definition}

It can be shown that an $\mathcal H^2$-matrix can be constructed with 
$\mathcal O(n k)$ computational and storage complexity. 
Consequently, the cost of a single matrix–vector product is also in 
$\mathcal O(n k)$, i.e., linear in the number of degrees of freedom.

Throughout this paper, we employ the \emph{GCA method} introduced by Börm and 
Christophersen in ~\cite{borm2016approximation}. 
This method is capable of constructing an $\mathcal H^2$-matrix approximation 
directly to a prescribed accuracy $\varepsilon_{\mathrm{ACA}} > 0$ and can thus 
be readily used to generate different accurate approximations of the integral operator 
defined in \eqref{eqn:integral_eq}.
Likewise other techniques can also be applied like interpolation
~\cite{borm2003introduction, borm2004H2, hackbusch2002data, Hackbusch2002, hackbusch2015hierarchical, borm2010efficient}
or adaptive cross-approximation ~\cite{ostrowski2006fast, bebendorf2008approximation, bebendorf2003adaptive, bebendorf2009recompression}.

Unless stated otherwise, we consider the complex-valued case, i.e., $\mathbb K = \mathbb C$, 
with the 3D Helmholtz kernel function
\[
g_\kappa(x,y) = \frac{\exp(\mathrm{i} \kappa \lVert x - y \rVert_2)}{4 \pi \lVert x - y \rVert_2},
\]
for a real-valued wave parameter $\kappa = 1$. 

Higher accuracy, corresponding to a smaller value of $\varepsilon_{\mathrm{ACA}}$, 
leads to a higher memory footprint of the constructed $\mathcal H^2$-matrix. 
This is because higher accuracy requires larger local ranks for the cluster bases, 
thereby increasing the overall memory consumption.

\section{\texorpdfstring{{$\mathcal H^2$}}{H^2}-matrix-vector multiplication}
\label{sec:h2_mat_vec_seq}

In 2002, Börm and Hackbusch ~\cite{hackbusch2002data} introduced the basic algorithm 
for computing a matrix–vector product of an $\mathcal H^2$-matrix. 
Its structure consists of three parts: the \emph{forward transformation}, 
the \emph{multiplication phase}, and the \emph{backward transformation}, 
as illustrated in \cref{alg:h2_mvm}.

\begin{algorithm}
\caption{$\mathcal H^2$-matrix-vector multiplication}
\label{alg:h2_mvm}
\begin{algorithmic}[1]
\Procedure{H2\_mvm}{$\alpha, G, x, y$}
\State{$\widehat y \gets 0 \, , \quad \widehat x \gets 0$}
\State{\Call{forward}{$\Root(\TJ), x, \widehat x$}}
\State{\Call{fast\_addeval}{$\alpha, G, \Root(\TI), \Root(\TJ), \widehat x, \widehat y, x, y$}}
\State{\Call{backward}{$\Root(\TI), y, \widehat y$}}
\EndProcedure
\end{algorithmic}
\end{algorithm}

The first and last steps operate solely on the cluster bases and, consequently, on the 
structure of the cluster trees. 
We aim for an efficient parallelization of these algorithms and first examine them 
before proceeding to the more computationally intensive multiplication phase in between. 
The forward transformation in its basic form is shown in \cref{alg:forward_backward} (left).

\begin{algorithm}
\caption{Forward transformation of an input vector $x$ into the transformed 
coefficient vector $(\widehat x_s)_{s \in \TJ}$, and backward transformation 
of a transformed coefficient vector $(\widehat y_t)_{t \in \TI}$ into the output vector $y$.}
\label{alg:forward_backward}
\begin{minipage}[t]{0.48\textwidth}
\begin{algorithmic}[1]
\Procedure{Forward}{$s, x, \widehat x$}
\If{$\sons(s) = \emptyset$}
\State{$\widehat x_s \gets W_s^* x$}
\Else
\For{$s' \in \sons(s)$}
\State{\Call{Forward}{$s', x|_{\hat s'}, \widehat x_{s'}$}}
\State{$\widehat x_s \gets \widehat x_s + E_{s'}^* \widehat x_{s'}$}
\EndFor
\EndIf
\EndProcedure
\end{algorithmic}
\end{minipage}
\hfill
\begin{minipage}[t]{0.48\textwidth}
\begin{algorithmic}[1]
\Procedure{Backward}{$t, y, \widehat y$}
\If{$\sons(t) = \emptyset$}
\State{$y \gets y + V_t \widehat y_t$}
\Else
\For{$t' \in \sons(t)$}
\State{$\widehat y_{t'} \gets \widehat y_{t'} + E_{t'} \widehat y_t$}
\State{\Call{Backward}{$t', y|_{\hat t'}, \widehat y_{t'}$}}
\EndFor
\EndIf
\EndProcedure
\end{algorithmic}
\end{minipage}
\end{algorithm}

During the forward transformation, data is propagated from the leaves up to the root 
of the cluster tree. 
Consequently, the recursive calls generating the vectors $\widehat x_{s'}$ can 
be executed in parallel. 
However, the subsequent multiplications with the transfer matrices $E_{s'}^*$ 
for $s' \in \sons(s)$ cannot be performed simultaneously, 
since this would result in a race condition when updating the vector $\widehat x_s$. 
Therefore, synchronization is required to ensure the correctness of the transformation.

In the case of the backward transformation (\Cref{alg:forward_backward}, right), 
the behavior differs slightly: during the recursive steps, we first update the vectors 
$\widehat y_{t'}$ and then propagate the information further down to the leaves. 
This makes things easier and allows the loop to be directly parallelized without any
need of synchronization on the same level of the hierarchy. 
The resulting parallel algorithms are shown in \cref{alg:par_forward_backward}.

\begin{algorithm}
\caption{Parallel versions of forward and backward transformations as a modification of \cref{alg:forward_backward}.}
\label{alg:par_forward_backward}
\begin{minipage}[t]{0.48\textwidth}
\begin{algorithmic}[1]
\Procedure{P\_Forward}{$s, x, \widehat x$}
\If{$\sons(s) = \emptyset$}
\State{$\widehat x_s \gets W_s^* x$}
\Else
\ParFor{$s' \in \sons(s)$}
\State{\Call{P\_Forward}{$s', x|_{\hat s'}, \widehat x_{s'}$}}
\EndParFor
\For{$s' \in \sons(s)$}
\State{$\widehat x_s \gets \widehat x_s + E_{s'}^* \widehat x_{s'}$}
\EndFor
\EndIf
\EndProcedure
\end{algorithmic}
\end{minipage}
\hfill
\begin{minipage}[t]{0.48\textwidth}
\begin{algorithmic}[1]
\Procedure{P\_Backward}{$t, y, \widehat y$}
\If{$\sons(t) = \emptyset$}
\State{$y \gets y + V_t \widehat y_t$}
\Else
\ParFor{$t' \in \sons(t)$}
\State{$\widehat y_{t'} \gets \widehat y_{t'} + E_{t'} \widehat y_t$}
\State{\Call{P\_Backward}{$t', y|_{\hat t'}, \widehat y_{t'}$}}
\EndParFor
\EndIf
\EndProcedure
\end{algorithmic}
\end{minipage}
\end{algorithm}

The major part of the runtime for a matrix–vector product arises from the 
multiplication phase. In its basic form, this procedure is shown in 
\cref{alg:h2_fast_addeval}.

\begin{algorithm}
\caption{Standard $\mathcal H^2$-matrix multiplication phase}
\label{alg:h2_fast_addeval}
\begin{algorithmic}[1]
\Procedure{fast\_addeval}{$\alpha, G, t, s, \widehat x_s, \widehat y_t, x, y$}
\If{$b = (t,s) \in \LIJ^+$}
\State{$\widehat y_t \gets \widehat y_t + \alpha\, S_b \widehat x_s$}
\ElsIf{$(t,s) \in \LIJ^-$}
\State{$y|_{\hat t} \gets y|_{\hat t} + \alpha\, G|_{\hat t \times \hat s} x|_{\hat s}$}
\Else
\ForAll{$(t', s') \in \sons(t,s)$}
\State{\Call{fast\_addeval}{$\alpha, G, t', s', \widehat x_{s'}, \widehat y_{t'}, x, y$}}
\EndFor
\EndIf
\EndProcedure
\end{algorithmic}
\end{algorithm}

Naively, all recursive calls within the loop could be executed in parallel. 
However, this would create a race condition among threads sharing the same row cluster $t$. 
Therefore, the algorithm must be modified such that matrix blocks associated with the 
same row cluster are executed sequentially on the same thread. 
This requires a preparation step that constructs lists of blocks for each row cluster 
$t \in \TI$, denoted by $\Clist_t$. 
In a second step, depending on whether a block is admissible or inadmissible, 
we perform matrix–vector multiplications of the corresponding sub-vectors.

\begin{figure}
\centering
\includegraphics[width=0.9\textwidth]{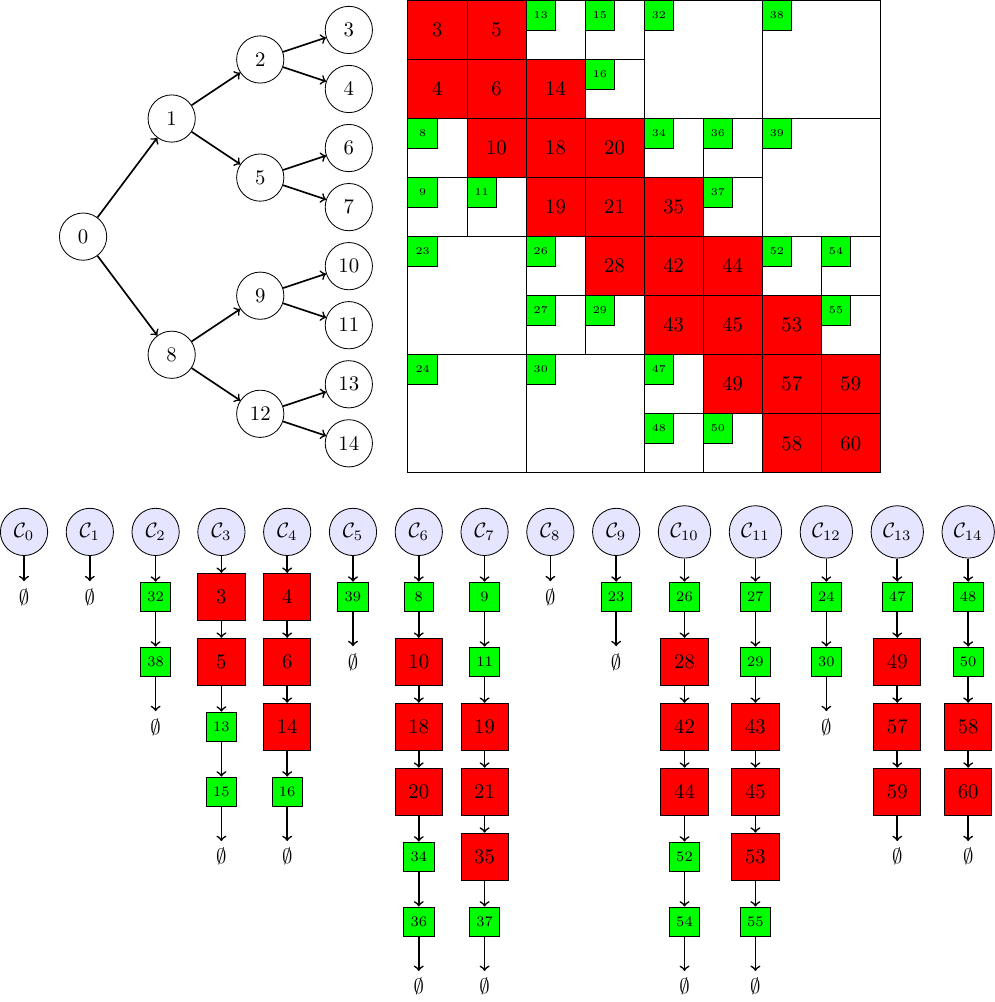}
\caption{Numbering of clusters and blocks, and their arrangement into the lists 
$\Clist_t$ for $t \in \TI$.
The cluster and block trees are numbered in \emph{pre-order}.
Green blocks are admissible, while red blocks are inadmissible.}
\label{img:list_mvm}
\end{figure}

The preparation step is described in \cref{alg:prep_mvm}, 
and the row-wise multiplication is carried out in \cref{alg:list_fast_mvm}.
An illustration of the preparation step can also be seen in
\cref{img:list_mvm}.

\begin{algorithm}
\caption{Preparation of the $\mathcal H^2$-matrix multiplication phase. 
         Sorting all leaf blocks into lists corresponding to their row cluster.
         Initially, $\Clist_t = \emptyset$ for all $t \in \TI$.}
\label{alg:prep_mvm}
\begin{algorithmic}[1]
\Procedure{prepare\_addeval}{$G, \Clist_t, s$}
\If{$\sons(G) \neq \emptyset$}
\ForAll{$(t', s') \in \sons(G)$}
\State{\Call{prepare\_addeval}{$G|_{\hat t' \times \hat s'}, \Clist_{t'}, s'$}}
\EndFor
\Else
\State{$\Clist_{t} \gets \Clist_{t} \cup (G, s)$}
\EndIf
\EndProcedure
\end{algorithmic}
\end{algorithm}

\begin{algorithm}
\caption{List-wise multiplication phase for each row cluster $t \in \TI$ separately.
All clusters can be handled in parallel and no synchronization between them is needed.}
\label{alg:list_fast_mvm}
\begin{algorithmic}[1]
\Procedure{list\_addeval}{$\alpha, \Clist_{t}, t, \widehat x, \widehat y, x, y$}
\If{$\sons(t) \neq \emptyset$}
\ParFor{$t' \in \sons(t)$}
\State{\Call{list\_addeval}{$\alpha, \Clist_{t'}, t', \widehat x, \widehat y, x, y$}}
\EndParFor
\Else
\While{$(G, s) \gets \operatorname{next}(\Clist_{t}) \neq \emptyset$}
\If{$b = (t,s) \in \LIJ^+$}
\State{$\widehat y_{t} \gets \widehat y_{t} + \alpha\, S_b \widehat x_{s}$}
\Else
\State{$y|_{\hat t} \gets y|_{\hat t} + \alpha\, G|_{\hat t \times \hat s} x|_{\hat s}$}
\EndIf
\EndWhile
\EndIf
\EndProcedure
\end{algorithmic}
\end{algorithm}

The creation of the auxiliary data structures $(\Clist_t)_{t \in \TI}$ requires some 
additional computation. Specifically, each leaf block must be added to a list, so 
the complexity is proportional to the number of leaf blocks, which can be estimated 
using the \emph{sparsity constant} $C_{\mathrm{sp}}$ as $\mathcal O(C_{\mathrm{sp}} \cdot |\TI|)$. 
Fortunately, this preparation step only needs to be performed once for a matrix 
and can be reused in subsequent matrix–vector products, for example, when the matrix 
is used in iterative Krylov subspace methods.

\subsection{Numerical examples}\label{sec:num_exp_H2_MVM_MM}

\begin{figure}
\includegraphics[width=\textwidth]{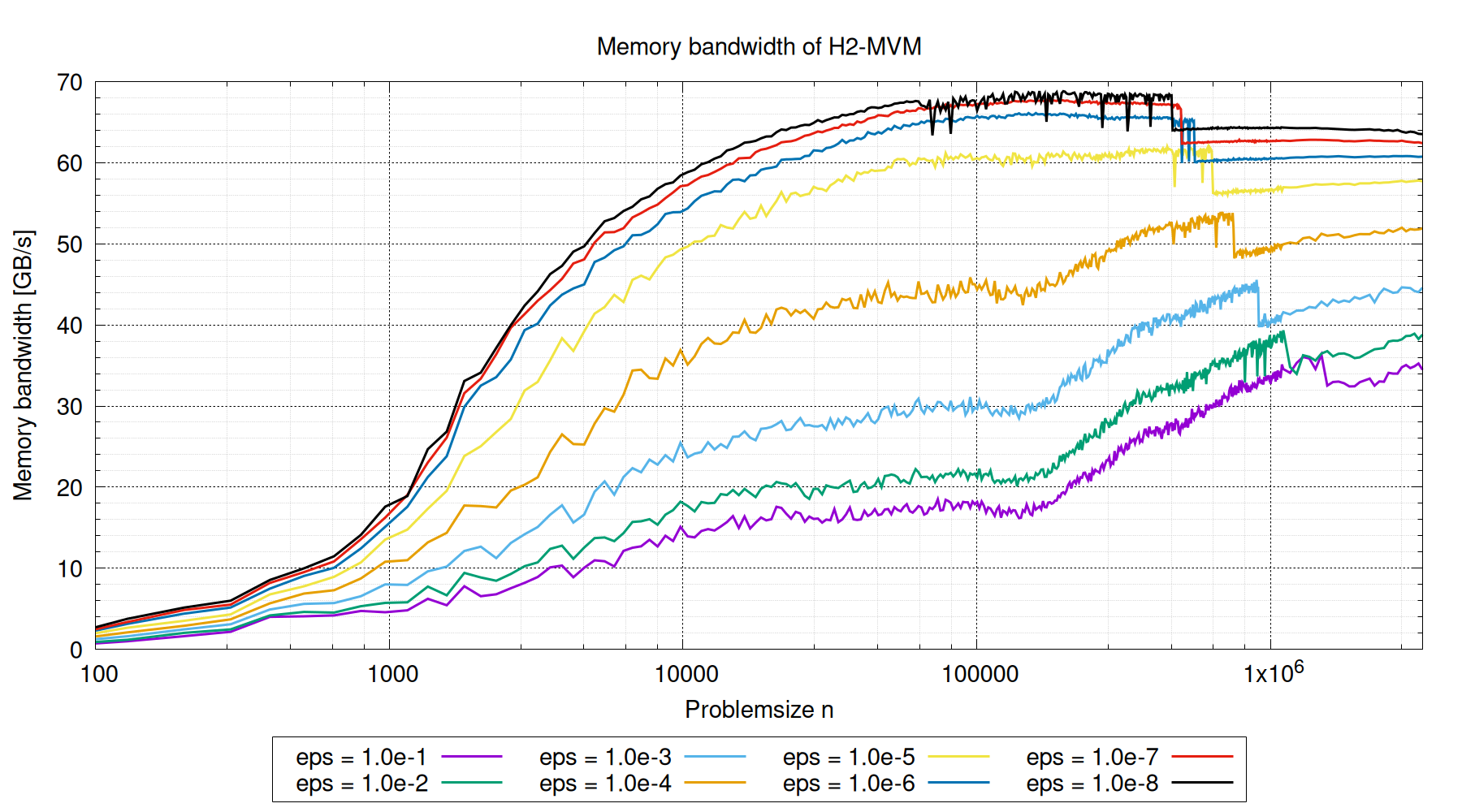}
\caption{Achievable bandwidth of parallel matrix–vector multiplication using 
\cref{alg:list_fast_mvm} with fast list evaluation.
48 threads are used, 24 physical + 24 logical cores of the CPU.}
\label{img:bandwidth_par_H2_MVM}
\end{figure}

Next, we discuss the achievable performance of the parallel approach. 
To this end, we present several numerical examples, the first of which is shown 
in \cref{img:bandwidth_par_H2_MVM}.

All experiments were performed on a single Intel Xeon 8160 Platinum processor, 
utilizing all 48 logical cores to maximize the total memory bandwidth. 
Benchmarking with the STREAM benchmark ~\cite{McCalpin2007} reports a total 
bandwidth of approximately 71~GB/s on this system. 
For the linear algebra subroutines we employ the intel MKL in version 2024.2.

As a test geometry for the integral equation~\eqref{eqn:integral_eq}, we use a sphere 
with varying numbers of triangles on its surface. 
\cref{img:bandwidth_par_H2_MVM} indicates that our algorithm can achieve 
this maximal bandwidth when enough threads are used to fully occupy the memory 
controller and sufficient data is available for loading. 
This occurs for a more accurate approximation of the integral operator, 
i.e., for a smaller $\varepsilon_{\mathrm{ACA}}$, such as $10^{-8}$. 

For other accuracies, and thus different memory footprints of the $\mathcal H^2$-matrix, 
the achieved bandwidth also increases with the problem size. 
Interestingly, in all test cases, the bandwidth temporarily drops at a certain point 
before increasing again. The explanation for this behavior remains unclear.

This Intel system features a 6-channel memory interface, so at most six threads 
are required in parallel to fully utilize the memory bandwidth. 
Since multiple read and write operations on the matrices must be performed simultaneously, 
it is reasonable to use a multiple of six threads to allow overlapping of these operations. 

\begin{figure}
\includegraphics[width=\textwidth]{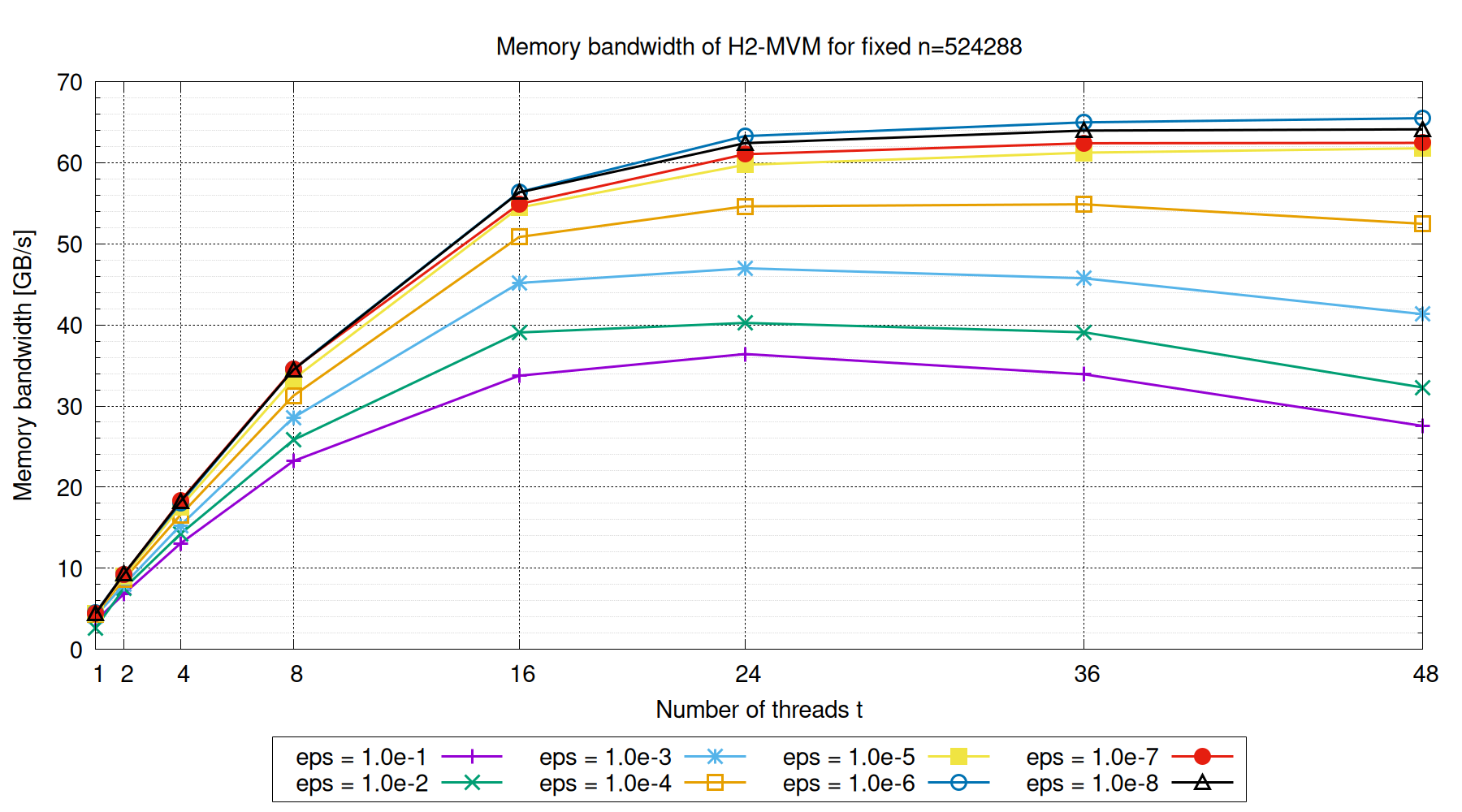}
\caption{Achievable bandwidth of parallel matrix–vector multiplication using 
\cref{alg:list_fast_mvm} with fast list evaluation for various numbers of threads.}
\label{img:bandwidth_thread_H2_MVM}
\end{figure}

\cref{img:bandwidth_thread_H2_MVM} shows the achievable transfer rates 
as a function of the number of threads. 
The bandwidth increases approximately linearly for a small number of threads, 
then gradually saturates until full utilization of 24 threads is reached. 
Using the logical hyper-threads of this CPU does not further increase bandwidth; 
for larger $\varepsilon_{\mathrm{ACA}}$ the bandwidth even decreases, 
while for smaller $\varepsilon_{\mathrm{ACA}}$ it remains roughly constant. 

In all cases, maximal bandwidth utilization is achieved for small $\varepsilon_{\mathrm{ACA}}$, 
corresponding to larger matrix blocks in the $\mathcal H^2$-matrix.

\section{\texorpdfstring{{$\mathcal H^2$}}{H^2}-matrix–matrix multiplication}

In a single matrix–vector multiplication, whether the structure of the matrix is dense, sparse, or an $\mathcal H^2$-matrix, 
each matrix entry is used only once, while the vector entries are read multiple times. 
This situation changes if we combine several vectors to form a dense matrix. 
Specifically, for $m \in \mathbb N$ and vectors $x^{(i)} \in \mathbb K^n$, $i = 1, \ldots , m$, 
let $\mathbf x \in \mathbb K^{n \times m}$ denote the resulting matrix. 
Similarly, by combining vectors $y^{(i)} \in \mathbb K^n$, $i = 1, \ldots , m$, 
we obtain a matrix $\mathbf y \in \mathbb K^{n \times m}$. 

Let $G \in \mathbb K^{n \times n}$ be an $\mathcal H^2$-matrix. 
Performing the operation
\[
\mathbf y \gets \mathbf y + \alpha\, G\, \mathbf x
\]
instead of the $m$ separate operations
\[
y^{(i)} \gets y^{(i)} + \alpha\, G\, x^{(i)}, \quad \text{for all } i = 1, \ldots , m,
\]
yields a matrix–matrix multiplication that can be implemented efficiently using
the BLAS level-3 routine \texttt{GEMM}.

Fortunately, multiplying an $\mathcal H^2$-matrix with a dense matrix is not fundamentally 
different from multiplying it with a single vector. 
We can reuse the preparation step from \cref{alg:prep_mvm} and adapt the parallel 
forward and backward transformations from \cref{alg:par_forward_backward} 
by applying them to submatrices $\mathbf x|_{\hat s \times m}$ and 
$\mathbf y|_{\hat t \times m}$, respectively. 
Similarly, the entities $\widehat x_s$ and $\widehat y_t$ are replaced by matrices with 
a small number of rows corresponding to the local rank of the cluster basis, 
i.e., $\mathbf{\widehat x_s} \in \mathbb K^{k_s \times m}$ and 
$\mathbf{\widehat y_t} \in \mathbb K^{k_t \times m}$. 

These modifications are sufficient in order to implement a fast and efficient 
$\mathcal H^2$-matrix–matrix multiplication. 
The performance of these algorithms is discussed in the next section.

\subsection{Numerical examples}

For the multiplication of an $\mathcal H^2$-matrix with a dense matrix having 
$m \in \mathbb N$ columns, the behavior differs slightly from the 
prior test case, where we had a single $\mathcal H^2$-matrix-vector product. 
Since the goal in this case is high data reuse, we cannot focus solely on the main memory bandwidth. 
At the same time, pure cache bandwidth is difficult to quantify. 
Therefore, we treat the multiplication of an $(n \times n)$ $\mathcal H^2$-matrix 
with a dense $n \times m$ matrix as $m$ separate matrix–vector multiplications with the same 
$\mathcal H^2$-matrix. 

Accordingly, the amount of data transferred is computed as the size of the $\mathcal H^2$-matrix 
multiplied by $m$, plus $m$ reads of the input vectors $x^{(i)}$ and $y^{(i)}$, 
and $m$ writes of the output vectors $y^{(i)}$. 

\begin{figure}
\begin{subfigure}{.5\textwidth}
\includegraphics[width=\textwidth]{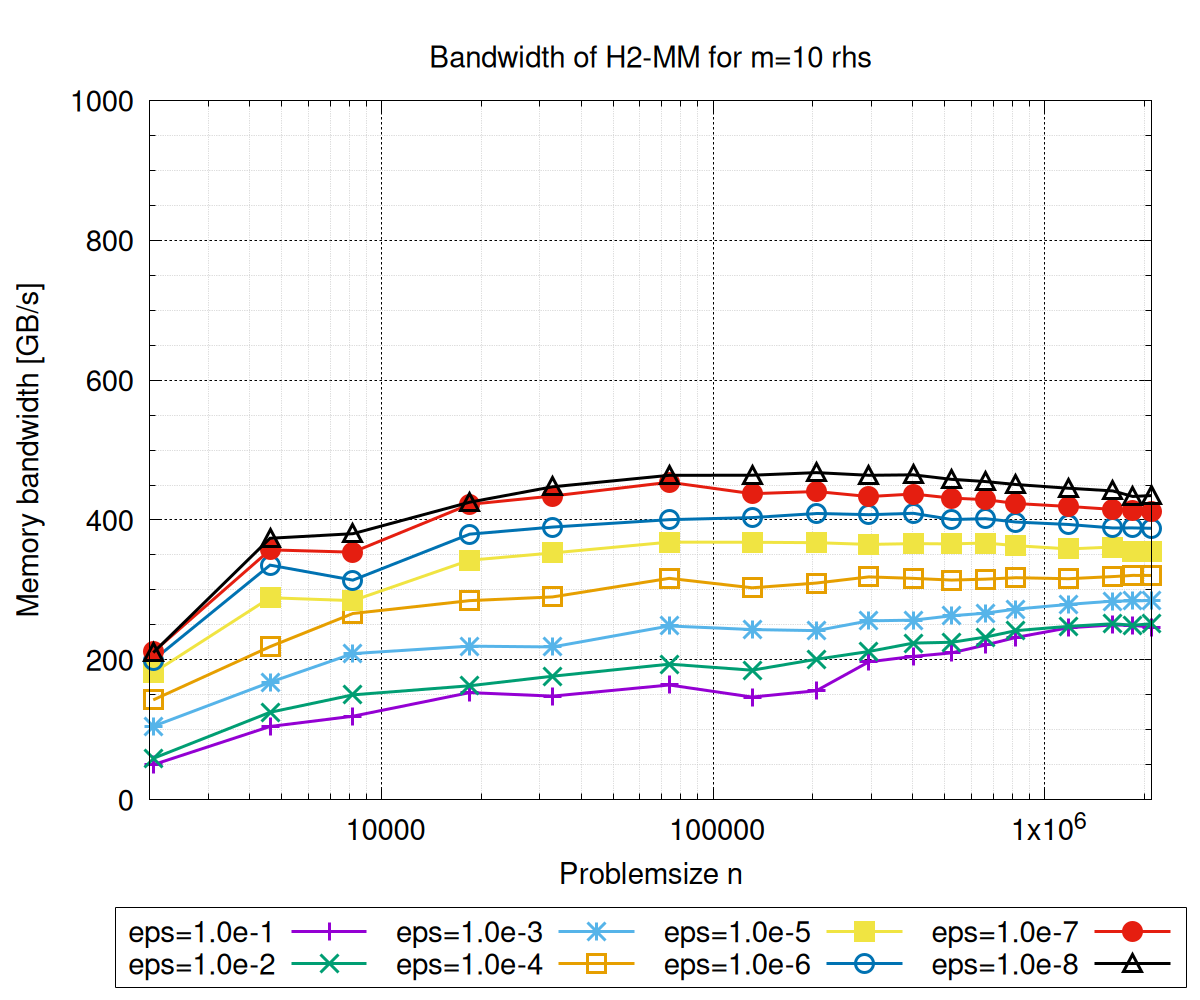}
\end{subfigure}%
\begin{subfigure}{.5\textwidth}
\includegraphics[width=\textwidth]{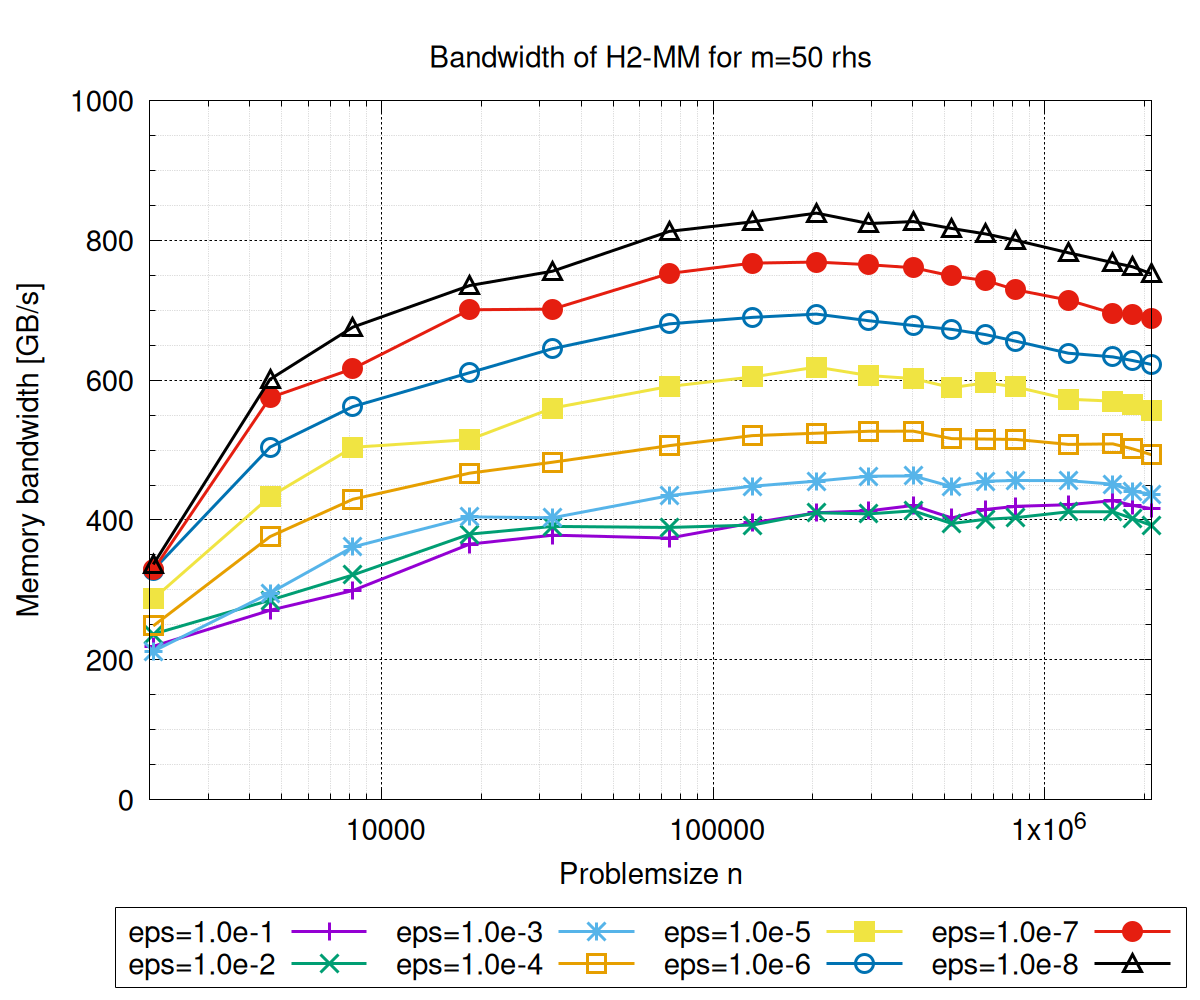}
\end{subfigure}%
\\
\begin{subfigure}{.5\textwidth}
\includegraphics[width=\textwidth]{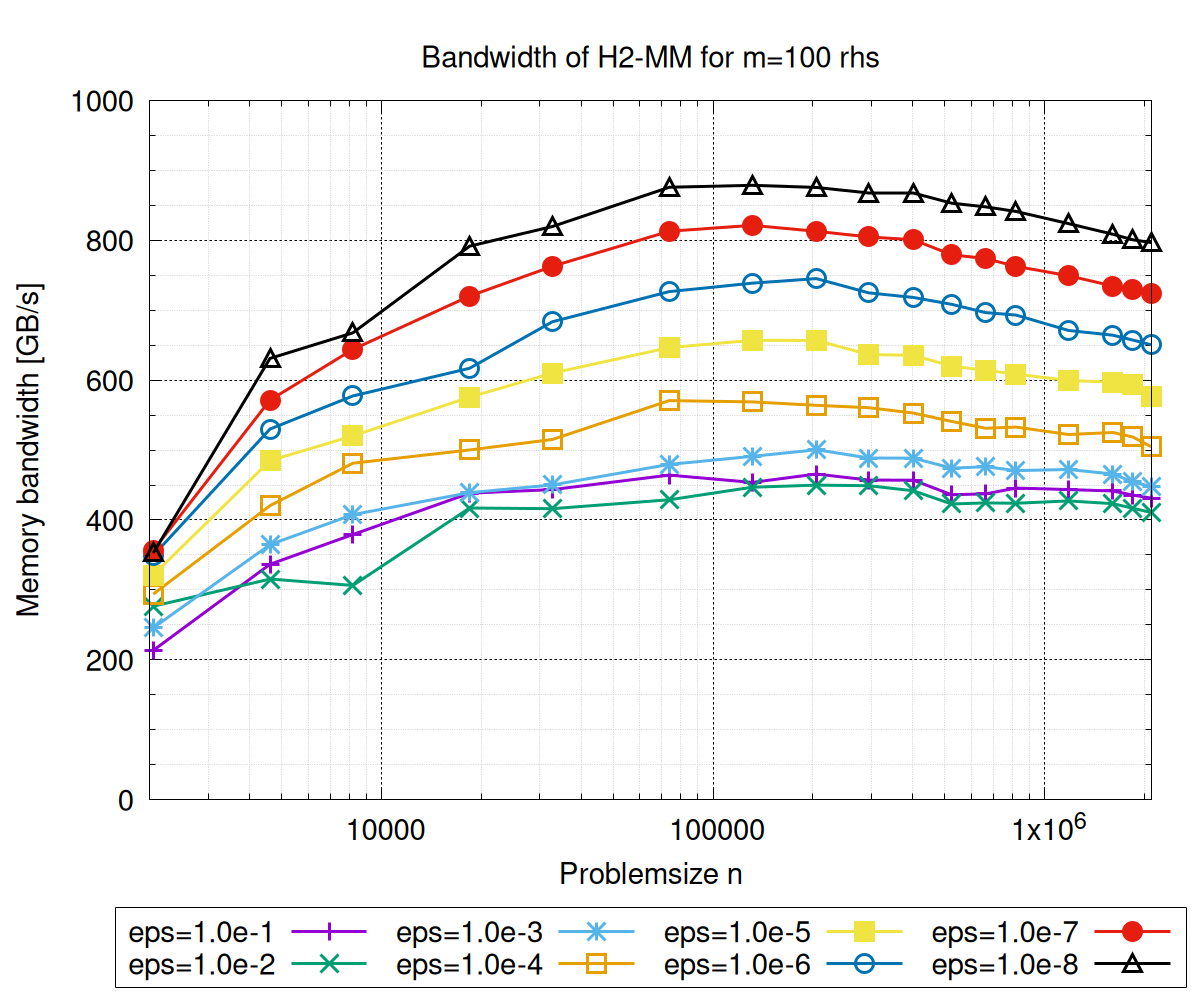}
\end{subfigure}%
\begin{subfigure}{.5\textwidth}
\includegraphics[width=\textwidth]{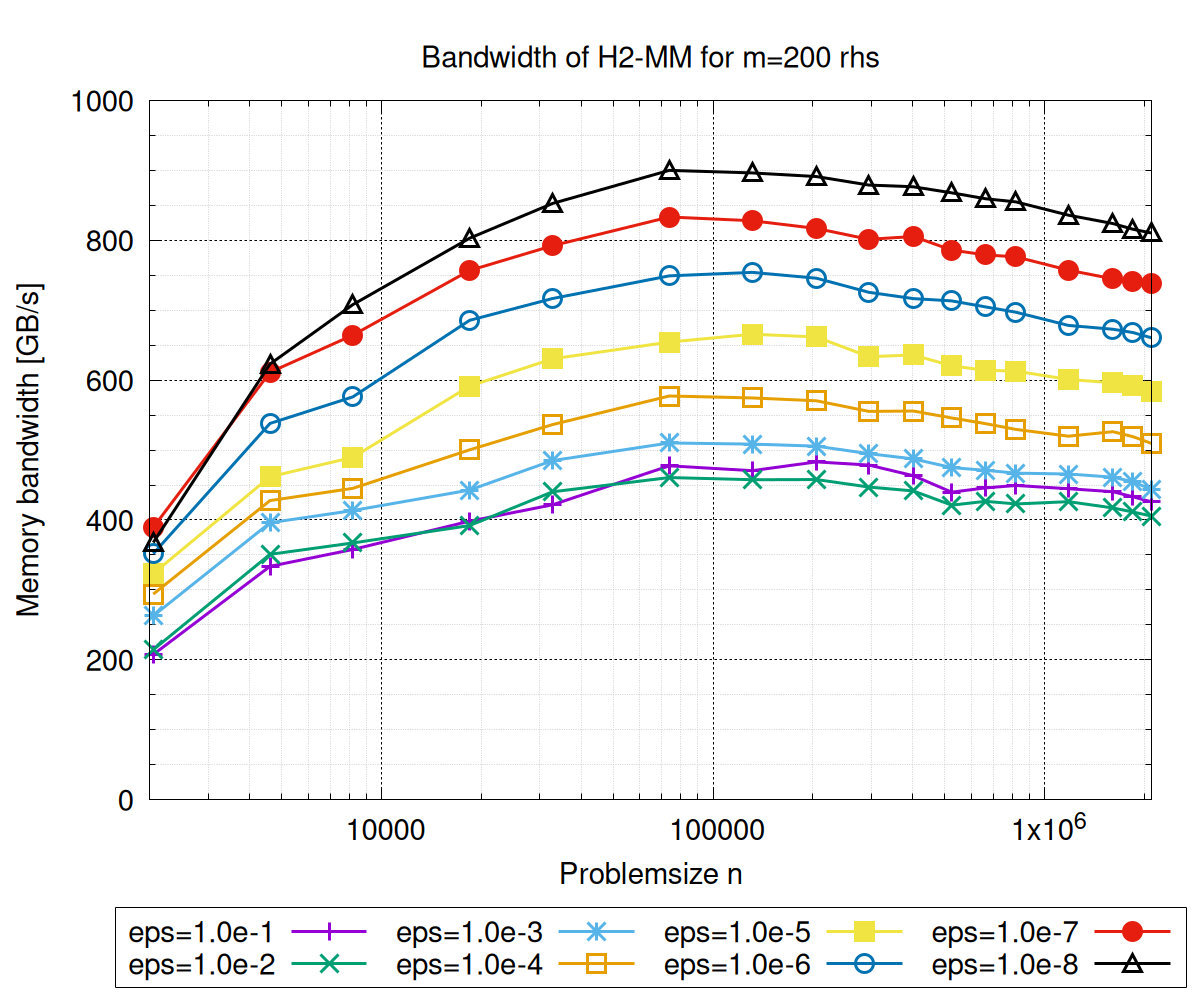}
\end{subfigure}%
\caption{Achievable bandwidth of parallel $\mathcal H^2$-matrix–matrix multiplication 
using a modified version of \cref{alg:list_fast_mvm} with fast list evaluation.
48 threads are used, 24 physical + 24 logical cores of the CPU.
The dense matrix has $m = 10$ columns(top left), $m = 50$(top right), 
$m = 100$(bottom left), and $m = 200$(bottom right).}
\label{img:bandwidth_H2_MM}
\end{figure}

As an initial overview, \cref{img:bandwidth_H2_MM} shows the achieved 
bandwidths for the matrix–matrix multiplication algorithm.

This figure exhibits effects similar to those observed for the matrix–vector product. 
The achievable bandwidth increases with the memory footprint of the $\mathcal H^2$-matrix, 
which corresponds to smaller values of $\varepsilon_{\mathrm{ACA}}$. 
Consequently, the measured transfer rates are higher by roughly an order of magnitude. 

As $m$ increases, the bandwidth also grows, reaching nearly the maximum for $m = 100$. 
Further increases in $m$ have little additional effect on the bandwidth. 
However, for matrix sizes larger than $100{,}000$, we observe a slight decrease in bandwidth 
utilization. This effect is present across all matrix approximations and for all values of $m$. 
It is strongly suspected that this behavior is caused by full utilization of the available 
main memory bandwidth, which prevents the cache from performing optimally in these cases.

\begin{figure}
\begin{subfigure}{.5\textwidth}
\includegraphics[width=\textwidth]{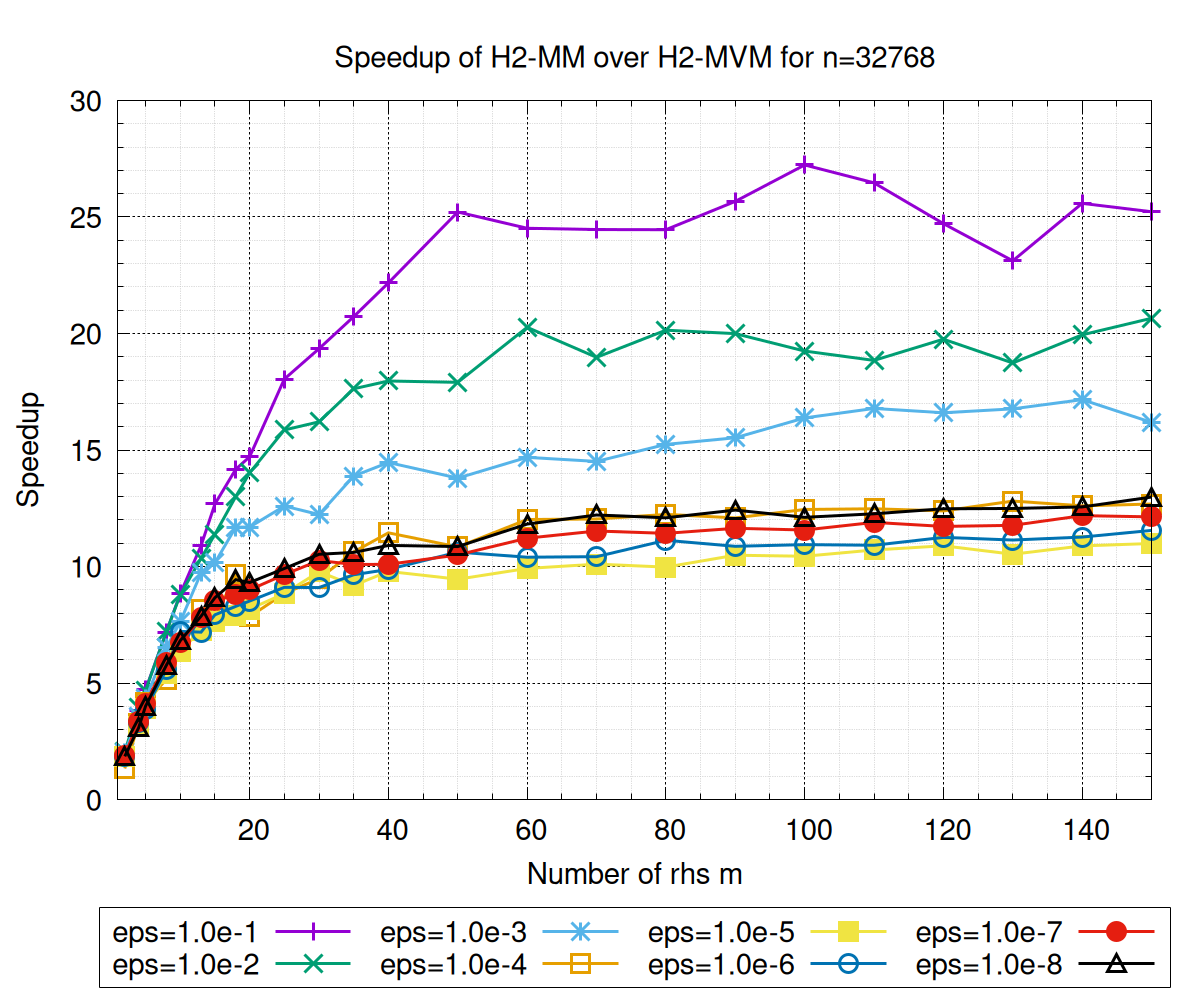}
\end{subfigure}%
\begin{subfigure}{.5\textwidth}
\includegraphics[width=\textwidth]{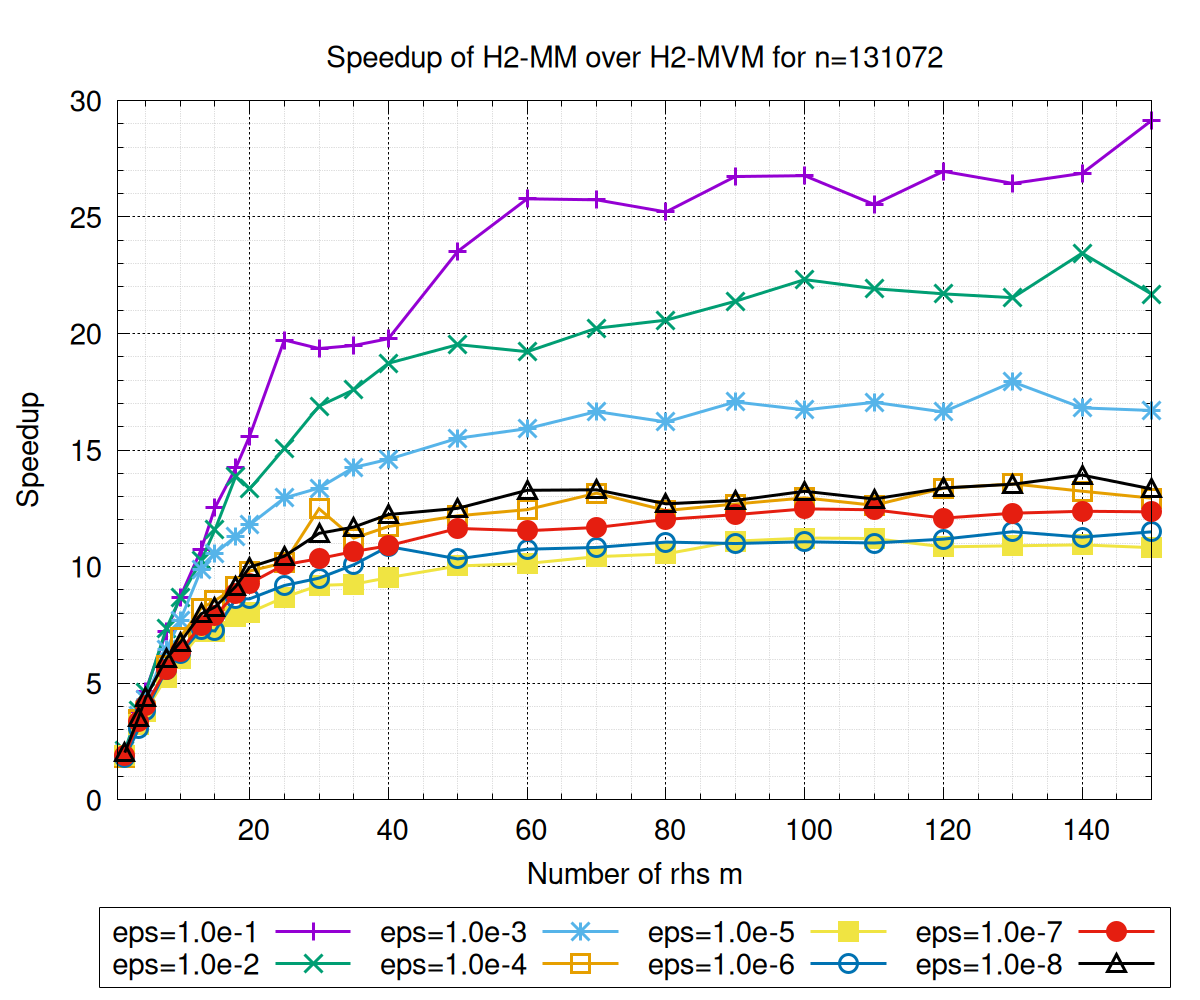}
\end{subfigure}%
\\
\begin{subfigure}{.5\textwidth}
\includegraphics[width=\textwidth]{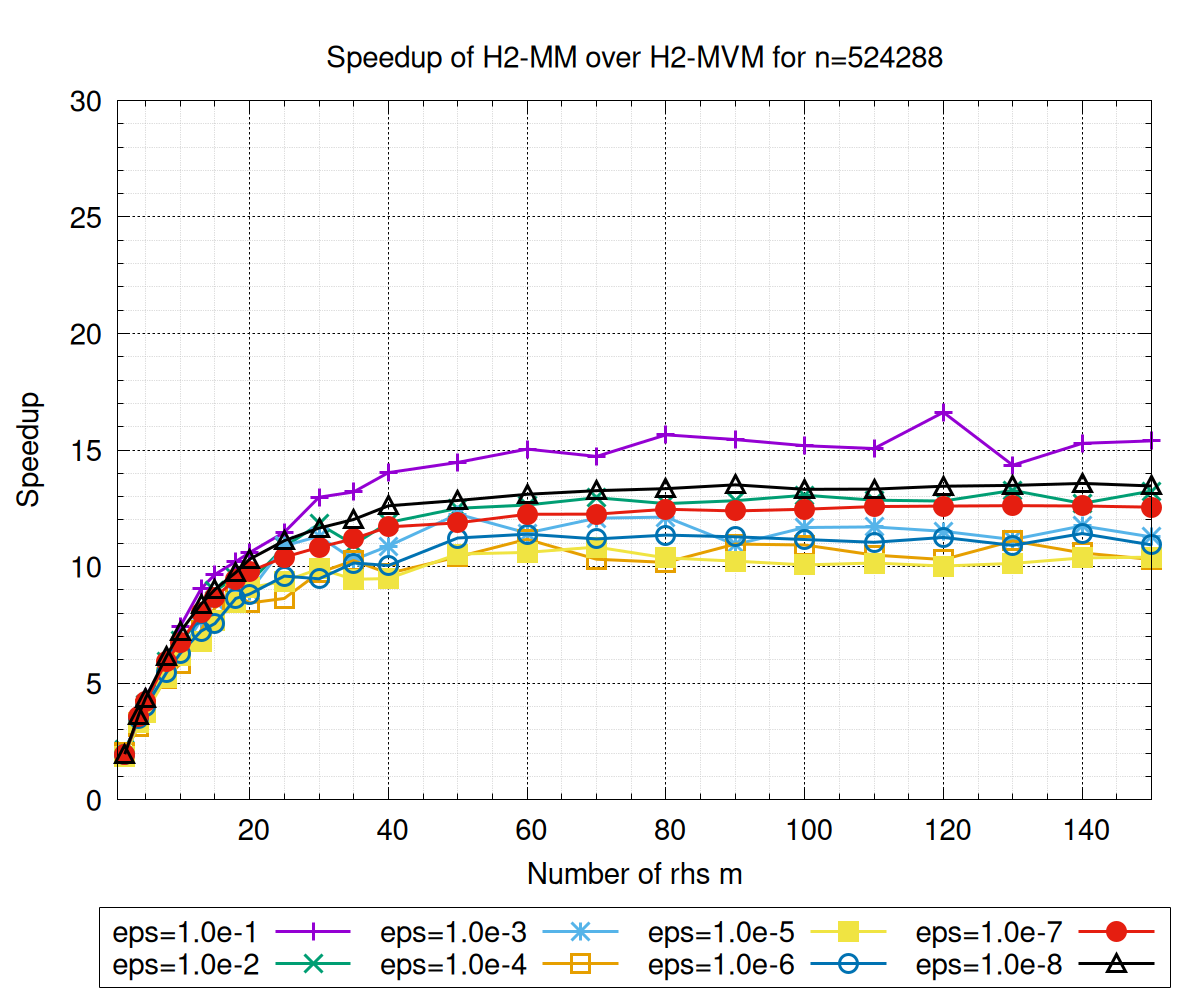}
\end{subfigure}%
\begin{subfigure}{.5\textwidth}
\includegraphics[width=\textwidth]{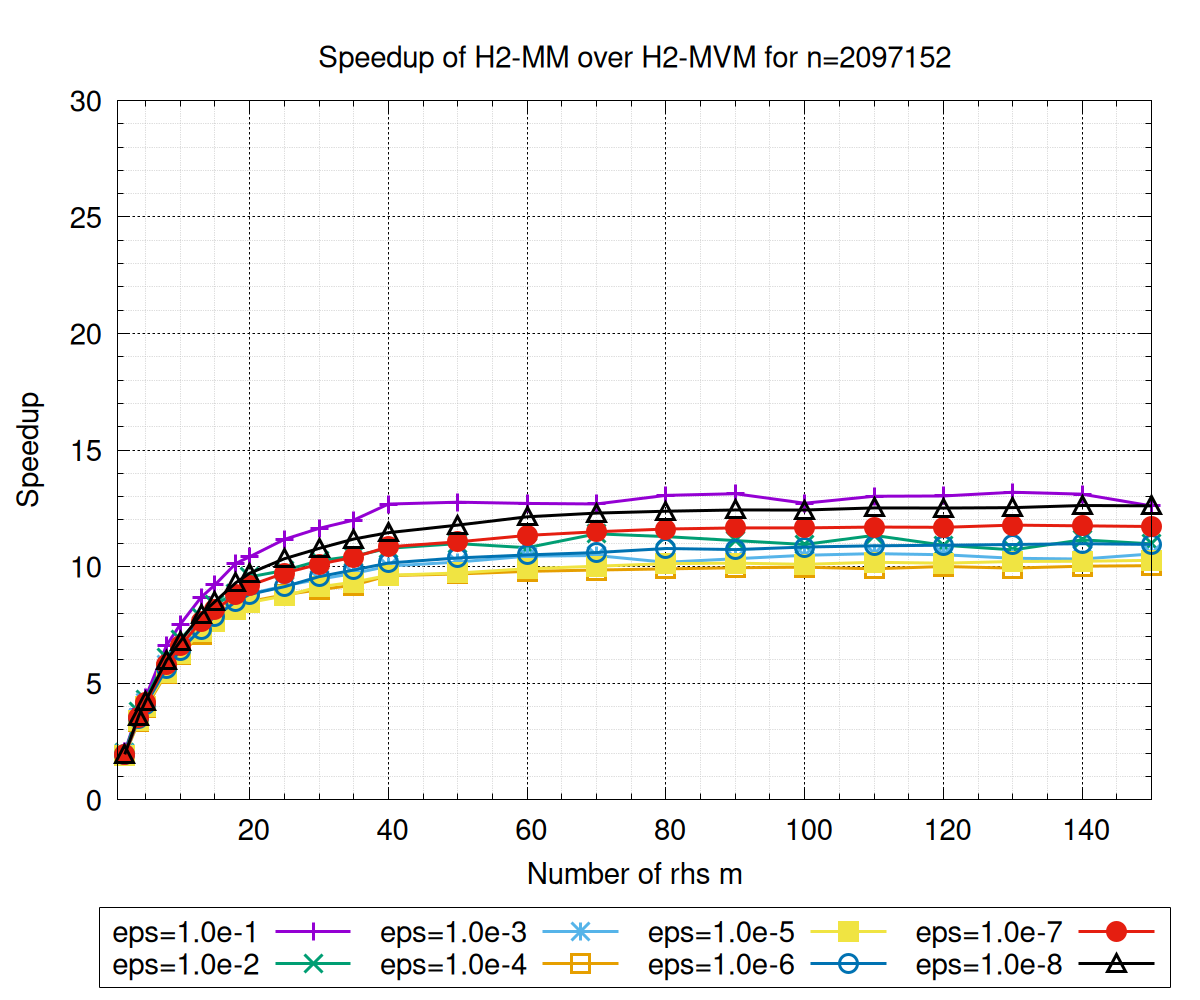}
\end{subfigure}%
\caption{Speedups of $\mathcal H^2$-matrix–matrix multiplication over 
$m$-fold $\mathcal H^2$-matrix–vector multiplication for different problem 
dimensions $n$ and varying numbers of right-hand sides $m$.}
\label{img:speedup_H2_MVM_vs_MM}
\end{figure}

Next, we'd like to draw the readers attention to \cref{img:speedup_H2_MVM_vs_MM}, 
which shows the possible speedups of the $\mathcal H^2$-matrix–matrix multiplication 
over the standard $m$-fold $\mathcal H^2$-matrix–vector product.

For small problem sizes, e.g., $n = 32{,}768$ (top left), more inaccurate 
matrix approximations, such as $\varepsilon_{\mathrm{ACA}} \in \{10^{-1}, 10^{-2}, 10^{-3}\}$, 
benefit more from intensive cache usage. 
This effect originates from the relatively poor performance of the matrix–vector product in these cases. 

\cref{img:speedup_H2_MVM_vs_MM} further indicates that for larger problems 
(bottom right), this is no longer the case, as the matrix–vector multiplication now 
utilizes the memory bandwidth more efficiently. 
Consequently, almost all accuracies exhibit similar speedups between 10 and 13 
for sufficiently many right-hand sides $m$.

Nevertheless, these achievable speedups of more than an order of magnitude are remarkable. 
Since matrix–vector products are frequently used as a key component in iterative solvers, 
such as Krylov subspace methods, our next goal is to incorporate this technique into 
such solvers.
This approach is discussed in the following section.
\newpage

\section{Block-Krylov Methods}
\label{sec:block_krylov}

We now investigate how standard Krylov methods can be adapted to efficiently handle 
the modified problem formulation of \eqref{eqn:matrix_system}. 
We aim to keep the $m$ linear systems independent of each other, i.e., 
each system evolves within its own Krylov subspace. 
Only shared operations, such as matrix–vector multiplication, are grouped and thus 
performed as matrix–matrix operations, namely the BLAS level-3 operation 
\texttt{GEMM}.

\subsection{Conjugate Gradient and Preconditioned CG}
\label{sec:cg}

We first revisit the well-known \emph{Conjugate Gradient (CG)} method, 
as outlined in \cref{alg:cg_alg}, where the system matrix $\mathbf A$ 
is symmetric and positive definite. 
In this setting, CG is one of the most efficient techniques for solving linear systems.

\begin{algorithm}
\caption{Conjugate gradient (CG) for a single right-hand side vector $b$}
\label{alg:cg_alg}
\begin{algorithmic}[1]
\State{Initialize $r \gets b - \mathbf A x, \quad p \gets r$}
\Comment{Initialize residual and search direction}
\While{$\lVert r \rVert_2 > \varepsilon_{\mathrm{slv}}$}
  \State{$a \gets \mathbf A p$}
  \Comment{Matrix–vector multiplication}
  \State{$\beta \gets \langle p, a \rangle_2, \quad 
         \alpha \gets \langle p, r \rangle_2 / \beta$}
  \Comment{Compute step size}
  \State{$x \gets x + \alpha p, \quad r \gets r - \alpha a$}
  \Comment{Update solution and residual}
  \State{$\gamma \gets \langle a, r \rangle_2 / \beta, \quad 
         p \gets r - \gamma p$}
  \Comment{Update search direction}
\EndWhile
\State{\Return $x$}
\end{algorithmic}
\end{algorithm}

On line~1, a single matrix–vector product is used to compute the initial residual, 
and within the while loop on line~3, another matrix–vector product is required. 
The key idea in a block Krylov method is to replace both operations with a 
matrix–matrix multiplication to improve cache utilization.

Transforming this algorithm into a fully blocked version also requires modifying 
the linear combinations of vectors. Since these are now columns of matrices, 
the operations must be performed column-wise. 
The same applies to the scalar products, whose results are now stored as vectors 
instead of single scalar values. 
For instance, the scalars $\alpha, \beta, \gamma \in \mathbb K$ are promoted to 
vectors $\boldsymbol \alpha, \boldsymbol \beta, \boldsymbol \gamma \in \mathbb K^m$. 
The resulting blocked version of \cref{alg:cg_alg} is shown in 
\cref{alg:cg_alg_block}.

\begin{algorithm}
\caption{Conjugate gradient method for $m$ right-hand sides stored in the matrix $\mathbf b$}
\label{alg:cg_alg_block}
\begin{algorithmic}[1]
\State{Initialize $\mathbf r \gets \mathbf b - \mathbf A \mathbf x, \quad 
       \mathbf p \gets \mathbf r$}
\Comment{Initialize residuals and search directions}
\While{not all systems have converged to accuracy $\varepsilon_{\mathrm{slv}}$}
  \State{$\mathbf a \gets \mathbf A \mathbf p$}
  \Comment{Matrix–matrix multiplication}
  \ParFor{$i = 1, \ldots , m$}
    \State{$\beta_i \gets \langle p^{(i)}, a^{(i)} \rangle_2, \quad 
           \alpha_i \gets \langle p^{(i)}, r^{(i)} \rangle_2 / \beta_i$}
    \Comment{Compute step size}
    \State{$x^{(i)} \gets x^{(i)} + \alpha_i p^{(i)}, \quad r^{(i)} \gets r^{(i)} - \alpha_i a^{(i)}$}
    \Comment{Update solution and residual}
    \State{$\gamma_i \gets \langle a^{(i)}, r^{(i)} \rangle_2 / \beta_i, \quad 
           p^{(i)} \gets r^{(i)} - \gamma_i p^{(i)}$}
    \Comment{Update search direction}
  \EndParFor
\EndWhile
\State{\Return $\mathbf x$}
\end{algorithmic}
\end{algorithm}

The condition of the while loop on line~2 can be implemented in various ways. 
In our approach, we require all linear systems to converge up to the prescribed accuracy
$\varepsilon_{\mathrm{slv}} > 0$. 
Consequently, the number of iterations is bounded by the maximum number of iterations 
among all $m$ systems. 
Alternative stopping criteria are also possible, for instance, terminating the procedure 
as soon as at least one system has converged. 
The choice heavily depends on the desired accuracy and the conditioning of the underlying 
linear systems.

Since all linear systems are treated independently, the for-loop on line~4 can be 
parallelized across the available computing cores to improve efficiency. 
The matrix–matrix multiplications benefit greatly from cache utilization and are 
parallelized using \cref{alg:list_fast_mvm} in its matrix version. 
This approach yields a highly optimized solver for multiple linear systems.

For practical applications, it is advantageous to include a preconditioner to accelerate 
the convergence of the systems. 
Hierarchical matrices provide a fast method to construct an approximate Cholesky 
factorization of the system matrix up to a prescribed accuracy in almost linear time. 
This produces a lower triangular matrix $\mathbf L \in \mathbb K^{n \times n}$ such that 
$\mathbf M = \mathbf L \mathbf L^* \approx \mathbf A$. 

Introducing an additional vector $q \in \mathbb K^n$ (or, in the case of block-CG, 
a matrix $\mathbf q \in \mathbb K^{n \times m}$ storing the preconditioned residual 
$\mathbf q = \mathbf M^{-1} \mathbf r$) leads directly to the preconditioned version 
of \cref{alg:cg_alg_block}, as shown in \cref{alg:pcg_alg_block}.

\begin{algorithm}
\caption{Preconditioned conjugate gradient method with preconditioner $\mathbf M$ 
for $m$ right-hand sides stored in the matrix $\mathbf b$}
\label{alg:pcg_alg_block}
\begin{algorithmic}[1]
\State{Let $\mathbf r \gets \mathbf b - \mathbf A \mathbf x, \quad 
       \mathbf q \gets \mathbf M^{-1} \mathbf r, \quad 
       \mathbf p \gets \mathbf r$}
\Comment{Initialize residuals, preconditioned residuals, and search directions}
\While{not all systems have converged to accuracy $\varepsilon_{\mathrm{slv}}$}
  \State{$\mathbf a \gets \mathbf A \mathbf p$}
  \Comment{Matrix–matrix multiplication}
  \ParFor{$i = 1, \ldots , m$}
    \State{$\beta_i \gets \langle p^{(i)}, a^{(i)} \rangle_2, 
           \quad \alpha_i \gets \langle p^{(i)}, r^{(i)} \rangle_2 / \beta_i$}
    \Comment{Compute step size}
    \State{$x^{(i)} \gets x^{(i)} + \alpha_i p^{(i)}, 
           \quad r^{(i)} \gets r^{(i)} - \alpha_i a^{(i)}$}
    \Comment{Update solution and residual}
  \EndParFor
  \State{$\mathbf q \gets \mathbf M^{-1} \mathbf r$}
  \Comment{Apply preconditioner to updated residuals}
  \ParFor{$i = 1, \ldots , m$}
    \State{$\gamma_i \gets \langle a^{(i)}, q^{(i)} \rangle_2 / \beta_i, 
           \quad p^{(i)} \gets q^{(i)} - \gamma_i p^{(i)}$}
    \Comment{Update search directions}
  \EndParFor
\EndWhile
\State{\Return $\mathbf x$}
\end{algorithmic}
\end{algorithm}

Again, on lines~1 and~8, the evaluation of the preconditioner $\mathbf M$ can be performed 
using BLAS level-3 operations and may also be parallelized to further improve the efficiency 
of the solver for multiple linear systems.

We would like to comment on some performance aspects of the Cholesky factorization: 
For a linear system $\mathbf M q = r$ with $\mathbf M = \mathbf L \mathbf L^*$ already 
factorized, one first solves $\mathbf L y = r$ and subsequently $\mathbf L^* q = y$. 
In the context of $\mathcal H$-arithmetic, this involves multiple small triangular 
matrix–vector solves (\texttt{TRSV}) and matrix–vector multiplications (\texttt{GEMV}). 
When processing $m$ right-hand sides simultaneously, these operations naturally transform 
into small triangular solves with multiple right-hand sides (\texttt{TRSM}) and matrix–matrix multiplications (\texttt{GEMM}), 
which are well-suited for high cache utilization.

It should be noted that the the triangular solves are sequential along the dependency chain, limiting parallelism; batching across right-hand sides restores parallel efficiency.
Since the $m$ right-hand sides are completely independent, we can regroup them: 
Assuming $p \in \mathbb N$ computing cores, we set $m_p := \frac{m}{p}$ and apply the 
preconditioner to chunks of $m_p$ right-hand sides at a time. 
This may slightly reduce cache efficiency, but ensures that all CPU cores receive a sufficient workload.

Before continuing to the next section, where we would like to discuss the same modification to 
the famous \emph{GMRES}-method, we at first present some numerical examples demonstrating the
performance of \cref{alg:pcg_alg_block} in practical applications, 
as well as a discussion of some limitations of this approach.

\subsubsection{Numerical Results}

\begin{figure}
\begin{subfigure}{.5\textwidth}
\includegraphics[width=\textwidth]{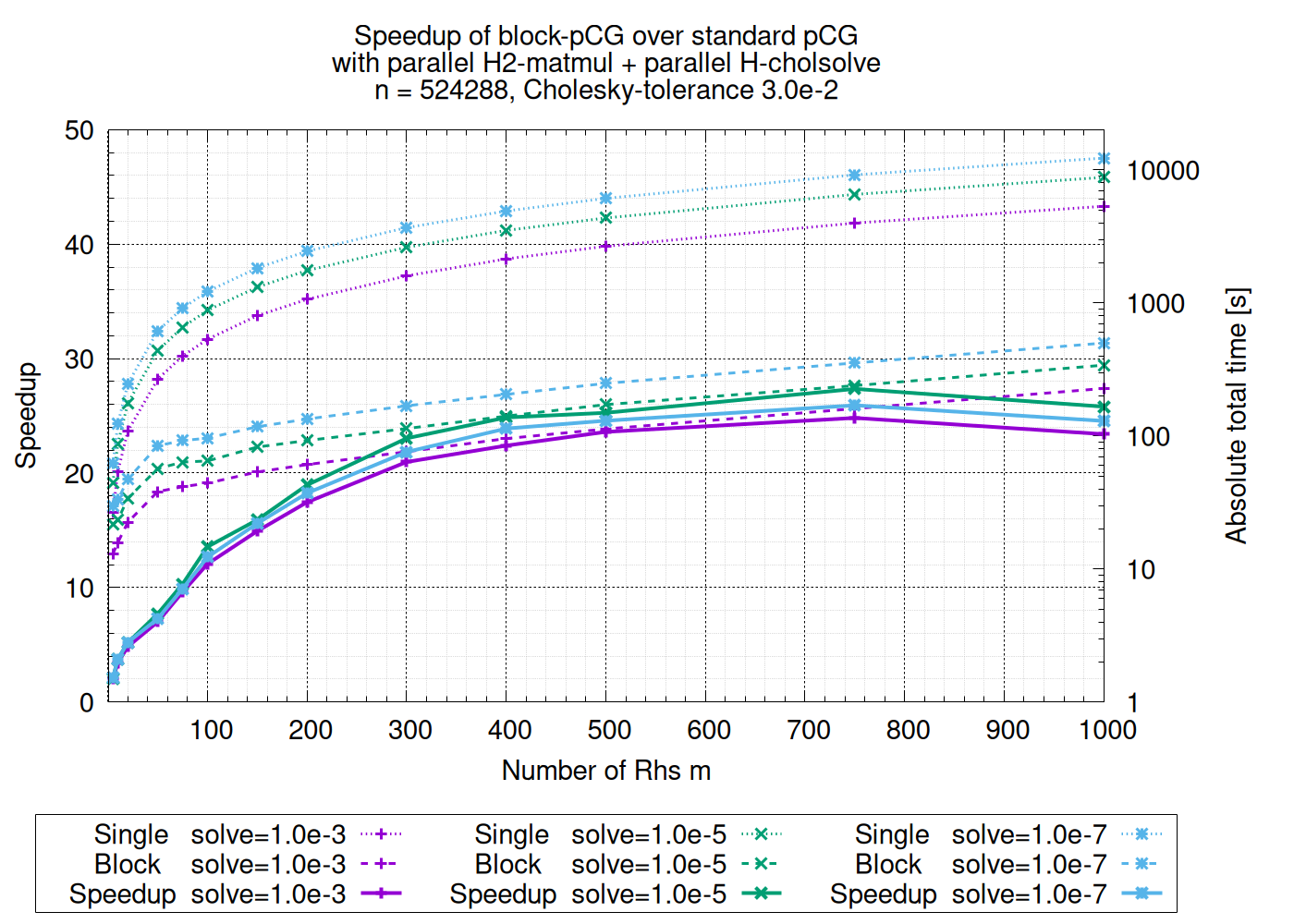}
\end{subfigure}%
\begin{subfigure}{.5\textwidth}
\includegraphics[width=\textwidth]{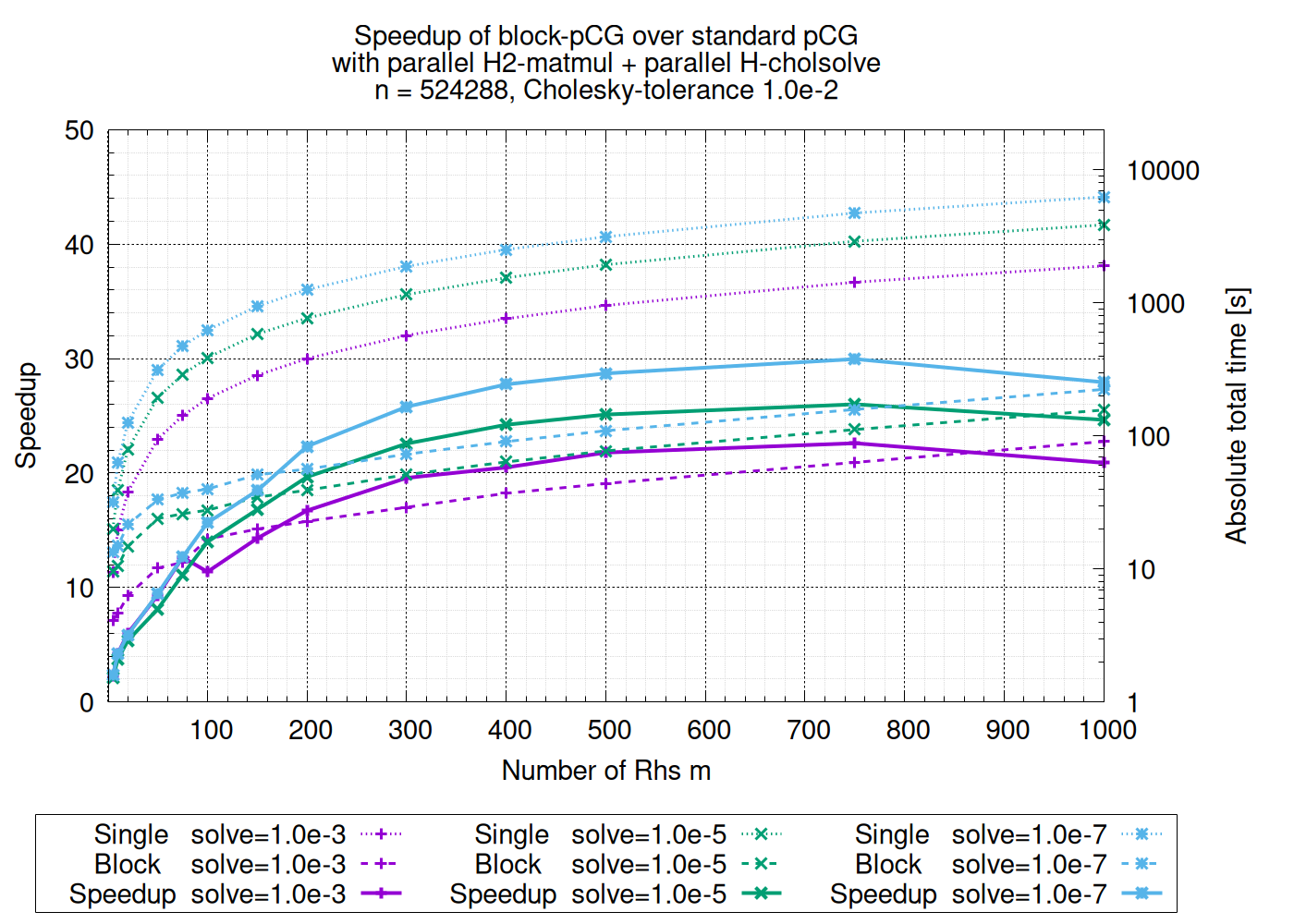}
\end{subfigure}%
\\
\begin{subfigure}{.5\textwidth}
\includegraphics[width=\textwidth]{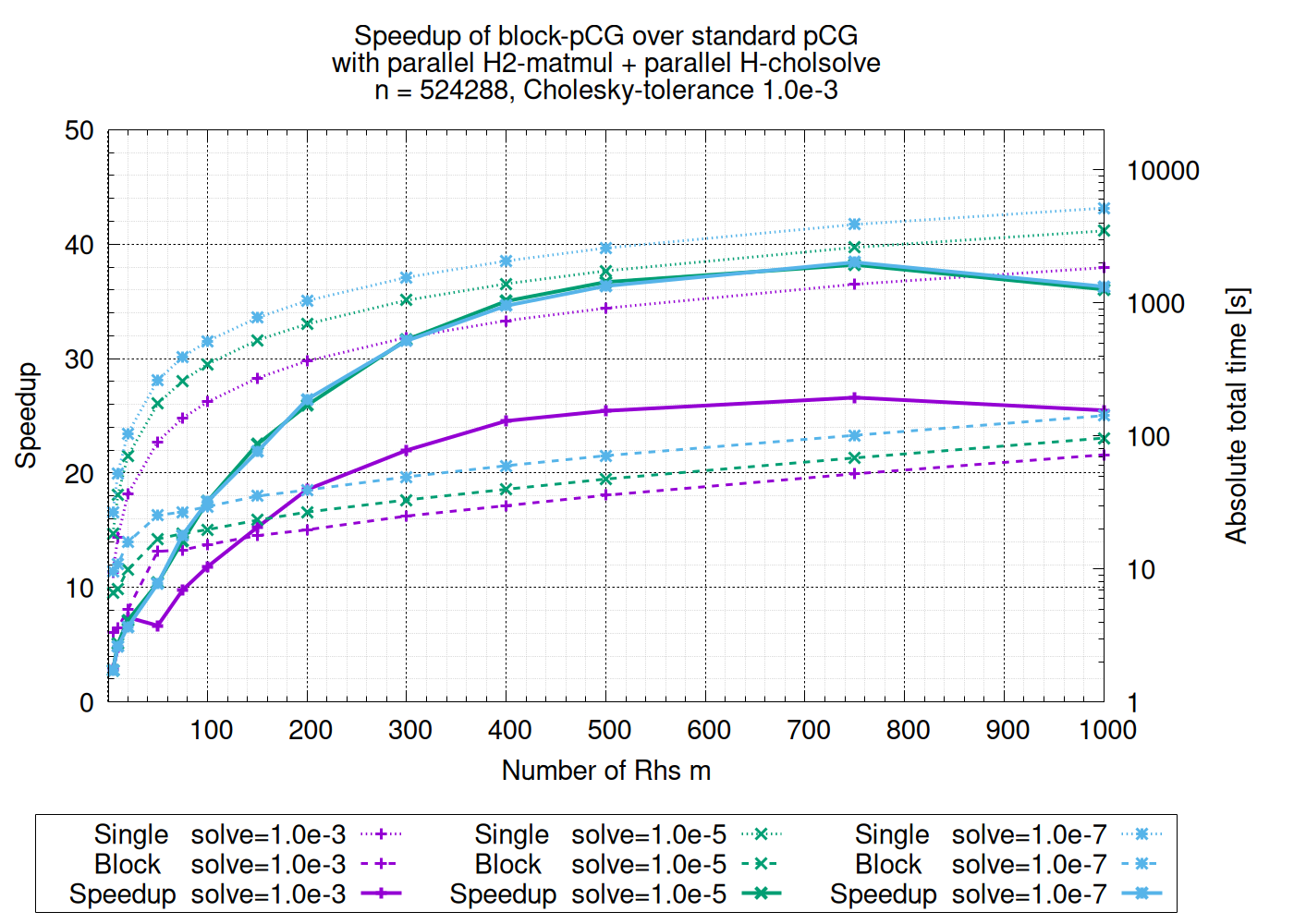}
\end{subfigure}%
\begin{subfigure}{.5\textwidth}
\includegraphics[width=\textwidth]{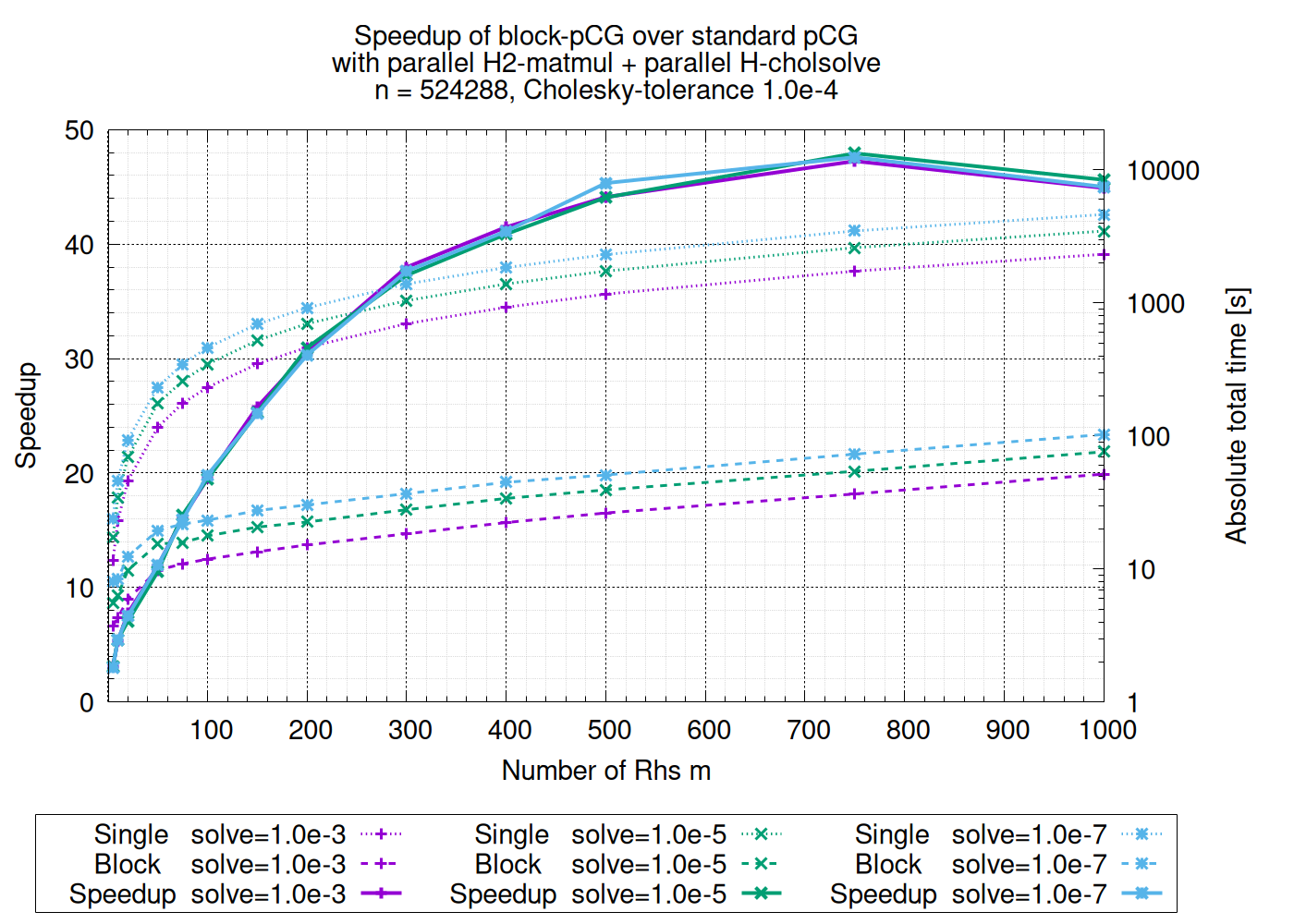}
\end{subfigure}%
\caption{Absolute runtimes for the standard pCG method and the block 
pCG method using parallel $\mathcal H^2$-matrix-vector multiplication and 
$\mathcal H^2$-matrix-matrix multiplication, respectively, along with the achievable 
speedups of the latter over the former. The accuracy of the Cholesky preconditioner 
$\varepsilon_{\mathrm{dcp}}$ and the relative solver tolerance $\varepsilon_{\mathrm{slv}}$ are varied.}
\label{img:speedup_H2_pCG}
\end{figure}

For this study, we employ the GCA method to construct an $\mathcal H^2$-matrix 
approximation of the single-layer operator for the real-valued Laplace equation. 
The operator is discretized on the unit sphere with an accuracy of $\varepsilon_{ACA} = 10^{-6}$, 
and the $\mathcal H$-Cholesky decomposition is used as a preconditioner for the system. 
The preconditioner itself is computed in $\mathcal H$-arithmetic up to accuracies 
$\varepsilon_{\mathrm{dcp}} \in \{10^{-1}, 10^{-2}, 10^{-3}, 10^{-4}\}$. 

To generate multiple linear systems, we create random solution vectors $x^{(i)}$, 
$i = 1, \ldots , m$, and compute the corresponding right-hand sides $b^{(i)} = \mathbf A x^{(i)}$. 
Thanks to the preconditioner, the $m$ systems exhibit similar convergence behavior, 
and the solution is achieved in an almost constant number of iterations using 
\cref{alg:pcg_alg_block}. 

The results, shown in \cref{img:speedup_H2_pCG}, indicate that similar performance 
is observed across different solver accuracies. 
Larger numbers of right-hand sides $m$ lead to improved cache utilization and more 
efficient parallelization of both the $\mathcal H^2$-matrix operations and the 
Cholesky preconditioner evaluation. 
In particular, the preconditioner benefits from an increased $m$, as discussed earlier, 
because a larger $m_p$ allows for better exploitation of the cache.

It can be observed that the speedup achieved in this setting is highest when a more accurate 
Cholesky preconditioner is used; in particular, $\varepsilon_{\mathrm{dcp}} = 10^{-4}$ yields the best performance. 
This is explained by the more efficient cache utilization in this case, which makes the 
application of the Cholesky factorization comparatively cheaper. 

Moreover, the speedups are almost identical across different solver tolerances $\varepsilon_{\mathrm{slv}}$,
except for the cases $\varepsilon_{\mathrm{dcp}} = 10^{-2}$ and
$\varepsilon_{\mathrm{dcp}} = 10^{-3}$. Here, the lowest solver accuracy appears
slightly less efficient than the others. This can be explained as follows:
since we require all $m$ linear systems to converge before the algorithm terminates,
some systems may need additional iterations even if they have individually reached
the prescribed accuracy, because others have not yet converged.
In such cases, the standard (non-blocked) algorithm may be slightly more advantageous.
This behavior can be observed in \cref{img:speedup_H2_pCG}.

It should also be noted that we were not able to construct a preconditioner
with an accuracy of $10^{-1}$.
The lowest achievable accuracy was $3 \cdot 10^{-2}$, as shown in
\cref{img:speedup_H2_pCG} (top left).

The maximum performance gain was observed for
$\varepsilon_{\mathrm{dcp}} = 10^{-4}$, where we achieved a speedup of
almost 50 compared to the standard approach of solving each linear system sequentially.

\begin{figure}
\begin{subfigure}{.5\textwidth}
\includegraphics[width=\textwidth]{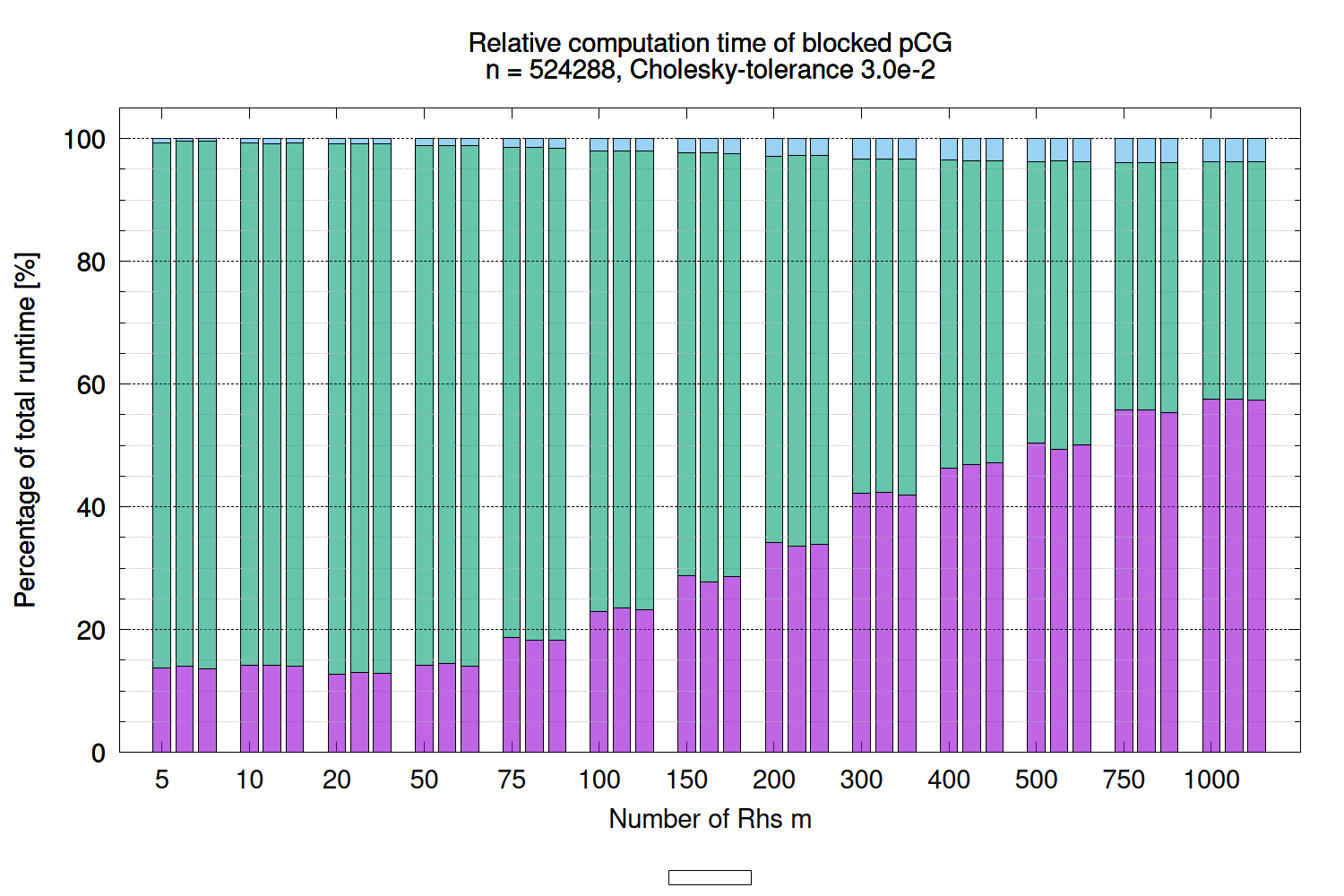}
\end{subfigure}%
\begin{subfigure}{.5\textwidth}
\includegraphics[width=\textwidth]{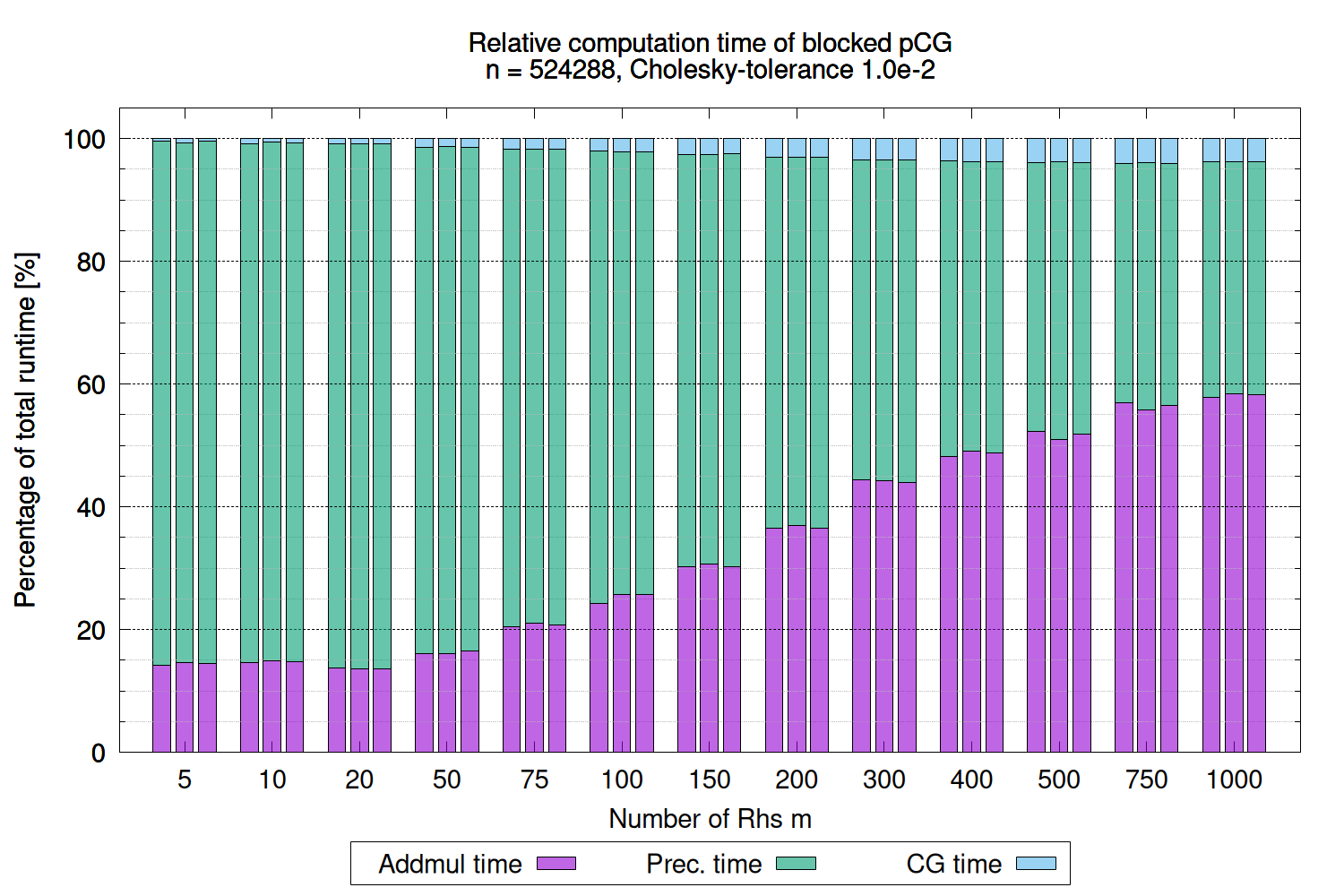}
\end{subfigure}%
\\
\begin{subfigure}{.5\textwidth}
\includegraphics[width=\textwidth]{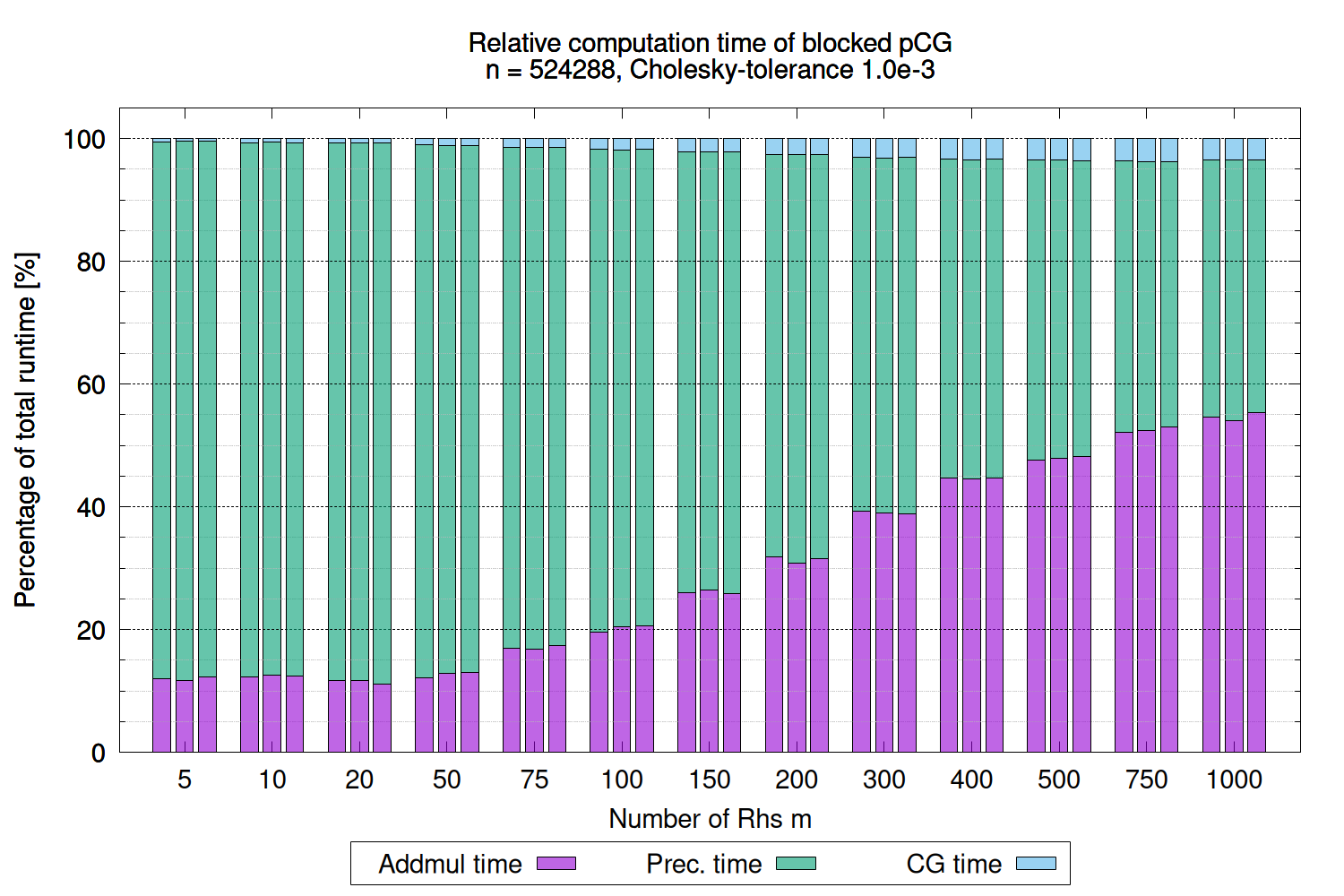}
\end{subfigure}%
\begin{subfigure}{.5\textwidth}
\includegraphics[width=\textwidth]{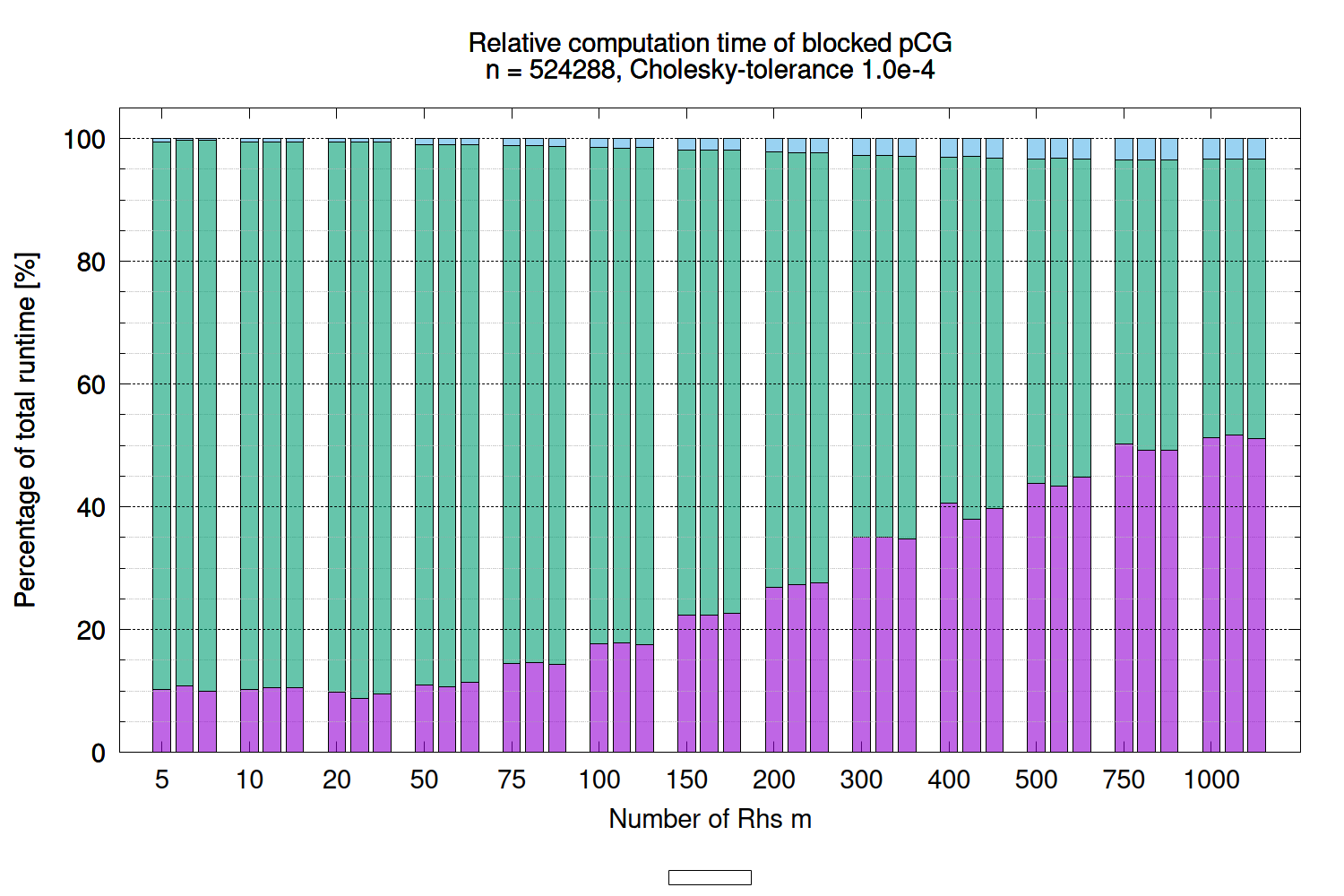}
\end{subfigure}%
\caption{Relative computing time of the solving process for \cref{alg:pcg_alg_block}. 
The accuracy of the Cholesky preconditioner $\varepsilon_{\mathrm{dcp}}$ and the solver tolerance 
$\varepsilon_{\mathrm{slv}}$ are varied. Within each group, the left bar corresponds to 
$\varepsilon_{\mathrm{slv}} = 10^{-3}$, the middle to $\varepsilon_{\mathrm{slv}} = 10^{-5}$, 
and the right to $\varepsilon_{\mathrm{slv}} = 10^{-7}$.}
\label{img:percentage_H2_pCG}
\end{figure}

Subsequently, we focus on the relative runtimes of the algorithm. 
We classify the $\mathcal H^2$-matrix-matrix multiplications in lines~1 and 3 of \cref{alg:pcg_alg_block} 
as ``addmul time''. Analogously, the evaluation of the preconditioner $\mathbf M$ in lines~1 and 8 
is considered as ``preconditioner time''. All remaining vector operations are grouped under ``CG time''. 
The percentage contributions of these three components to the total runtime are shown in \cref{img:percentage_H2_pCG}.

For a small number of right-hand sides $m$, applying the preconditioner accounts for roughly 85\% 
of the total runtime, while the $\mathcal H^2$-matrix-matrix multiplication consumes the remaining 15\%. 
As $m$ increases, the computation time for the matrix-matrix multiplication becomes dominant, 
and the ``CG time'' rises slightly to about 5\%. 
A higher $m$ also enables better parallelization of the preconditioner, thereby reducing its 
relative contribution to the total computation time.
All operation in the ``CG time'' section only work with BLAS level-1 function
and thus do not benefit from cache usage.

\begin{remark}{Memory footprint}{}
\label{rem:memory_block_pcg}
For a single linear system, the CG or pCG method requires storing a small number of vectors: 
$r$ for the residual, $p$ for the search direction, $a$ for the matrix–vector product, 
$q$ for the preconditioned residual, $b$ for the right-hand side, and $x$ for the solution vector. 
In total, this amounts to a memory usage of $6n$ floating point numbers. 
Additionally, a few scalar values such as $\alpha$, $\beta$, and $\gamma$ are needed, 
but their contribution can be neglected.  

The memory footprint of an $\mathcal H^2$-matrix is $\mathcal O(n k)$, 
where the rank $k$ depends on the approximation accuracy. 
Since typically $k \gg 6$, the $\mathcal H^2$-matrix itself dominates the memory consumption of the solver.  

The situation changes significantly when solving $m$ linear systems simultaneously. 
In this case, the vector quantities are promoted to matrices of size $n \times m$, 
requiring at least $\mathcal O(6nm)$ units of storage. 
Again, the vectors $\alpha$, $\beta$, and $\gamma$ can be neglected.  
As a consequence, for large values of $m$, the memory required by these auxiliary matrices 
can exceed that of the system matrix. 
Thus, $m$ must be chosen carefully. 
If $m$ becomes very large, it is advisable to process the systems in chunks of at most 
$m_{\mathrm{chunk}} < m$ right-hand sides. 
\end{remark}

\subsection{GMRES and pGMRES}
\label{sec:gmres}

The same concept extends naturally to other Krylov subspace methods. 
Here, we illustrate it for the \emph{Generalized Minimal Residual (GMRES)} method. 
Like CG, GMRES requires a single matrix-vector product per iteration along with a few 
vector operations and scalar products. 
A notable difference is the \emph{Arnoldi} process, which is invoked in each iteration 
to generate a new orthogonal search direction for the subsequent step. 

For efficiency, GMRES is typically limited to at most $\ell \in \mathbb N$ iterations. 
If convergence is not achieved within $\ell$ steps, the method can be restarted with an 
improved initial guess. This strategy is motivated both by efficient memory management 
and the fact that later iterations become increasingly expensive. 
The basic GMRES algorithm is presented in \cref{alg:gmres_alg}.

\begin{algorithm}
\caption{Generalized minimal residual method for a single right-hand side $b$
with an initial guess $x$. The updated solution is stored in $x$.
A maximum of $\ell$ iterations is performed unless the residual error falls below
$\varepsilon_{\mathrm{slv}}$.}
\label{alg:gmres_alg}
\begin{algorithmic}[1]
\State{Let $k = 0$}
\State{Let $r \gets b - \mathbf{A} x$} \Comment{initial, unnormalized residual}
\State{Let $\beta \gets \lVert r \rVert_2$} \Comment{initial residual norm}
\State{Let $\hat r \gets \beta e_1 \in \mathbb K^\ell$} \Comment{transformed initial residual}
\State{Let $q^{(k)} \gets r / \beta, \quad \mathbf V^{(k)} \gets \begin{pmatrix}q^{(k)}\end{pmatrix}$} 
\Comment{first basis vector}
\State{Let $\mathbf H \gets 0 \in \mathbb K^{(\ell + 1) \times \ell}$}
\Comment{initialize Hessenberg matrix with zeros}
\While{$k < \ell$ and $\lvert \hat r_k \rvert > \varepsilon_{\mathrm{slv}}$}
\State{$a^{(k+1)} \gets \mathbf A q^{(k)}$} \Comment{compute new Krylov vector}
\State{$\mathbf H|_{k \times \{k\}} \gets \mathbf{(V^{(k)})^*} a^{(k+1)}$}
\Comment{expand Hessenberg matrix}
\State{$a^{(k+1)} \gets a^{(k+1)} - \mathbf V^{(k)} \mathbf H|_{k \times \{k\}}$}
\Comment{orthogonalize new direction}
\State{Let $h_{k+1, k} \gets \lVert a^{(k+1)} \rVert_2$}
\Comment{finalize new column of Hessenberg matrix}
\State{Let $q^{(k+1)} \gets a^{(k+1)} / h_{k+1, k}$}
\Comment{normalize new direction}
\State{Let $\mathbf V^{(k+1)} \gets \begin{pmatrix}\mathbf V^{(k)} & q^{(k+1)}\end{pmatrix}$}
\Comment{update orthonormal basis}
\State{Construct Givens rotation $\mathbf G_k$ to eliminate $h_{k+1,k}$}
\State{$\mathbf H \gets \mathbf G_k \mathbf H$}
\Comment{apply $\mathbf G_k$ to $\mathbf H$}
\State{$\hat r \gets \mathbf G_k \hat r$}
\Comment{apply $\mathbf G_k$ to $\hat r$}
\State{$k \gets k + 1$}
\EndWhile
\State{Solve $\mathbf{R^{(k)}} y = \hat r$}
\Comment{$\mathbf{R^{(k)}}$: upper triangular part of $\mathbf H|_{k \times k}$ after Givens rotations}
\State{Update $x \gets x + \mathbf V^{(k)} y$}
\Comment{update $x$ with least-squares solution}
\State{\Return $x$}
\end{algorithmic}
\end{algorithm}

Analogously to the CG method, we now consider a preconditioned version of GMRES. 
In addition, we extend the algorithm to simultaneously handle $m$ right-hand sides, 
thereby exploiting BLAS level-3 operations for both the multiplication with the system matrix 
$\mathbf A$ and the application of the preconditioner $\mathbf M$. 
The resulting block pGMRES algorithm is presented in 
\cref{alg:pgmres_alg}.

\begin{algorithm}
\caption{Preconditioned generalized minimal residual method for a single right-hand side $b$
with an initial guess $x$. The updated solution is stored in $x$.
A maximum of $\ell$ iterations is performed unless the residual error falls below
$\varepsilon_{\mathrm{slv}}$.}
\label{alg:pgmres_alg}
\begin{algorithmic}[1]
\State{Let $k = 0$}
\State{Let $\mathbf r \gets \mathbf M^{-1} (\mathbf b - \mathbf{A} \mathbf x)$}
\Comment{initial, unnormalized residual matrix (column-wise)}
\ParFor{$i = 1, \ldots , m$}
\State{Let $\beta_i \gets \lVert {r_{*,i}} \rVert_2$}
\Comment{initial norm of the $i$-th residual vector}
\State{Let $\hat r^{(i)} \gets \beta_i e_1 \in \mathbb K^\ell$}
\Comment{transformed initial residual for system $i$}
\State{Let $q^{(k, i)} \gets r^{(i)} / \beta_i, \quad \mathbf V^{(k, i)} \gets \begin{pmatrix}q^{(k, i)}\end{pmatrix}$} 
\Comment{first basis vector for system $i$}
\State{Let $\mathbf H^{(i)} \gets 0 \in \mathbb K^{(\ell + 1) \times \ell}$}
\Comment{initialize Hessenberg matrix for system $i$ with zeros}
\EndParFor
\While{$k < \ell$ and not all systems have converged to accuracy $\varepsilon_{\mathrm{slv}}$}
\State{$\mathbf a^{(k+1)} \gets \mathbf M^{-1} \mathbf A \mathbf q^{(k)}$} 
\Comment{compute new Krylov vectors}
\ParFor{$i = 1, \ldots , m$}
\State{$\mathbf H^{(i)}|_{k \times \{k\}} \gets \mathbf{(V^{(k, i)})^*} a^{(k+1, i)}$}
\Comment{expand Hessenberg matrix for system $i$}
\State{$a^{(k+1, i)} \gets a^{(k+1, i)} - \mathbf V^{(k, i)} \mathbf H^{(i)}|_{k \times \{k\}}$}
\Comment{orthogonalize new direction}
\State{Let $h^{(i)}_{k+1, k} \gets \lVert a^{(k+1, i)} \rVert_2$}
\Comment{finalize new column of Hessenberg matrix}
\State{Let $q^{(k+1, i)} \gets a^{(k+1, i)} / h^{(i)}_{k+1, k}$}
\Comment{normalize new direction}
\State{Let $\mathbf V^{(k+1, i)} \gets \begin{pmatrix}\mathbf V^{(k, i)} & q^{(k+1, i)}\end{pmatrix}$}
\Comment{update orthonormal basis for system $i$}
\State{Construct Givens rotation $\mathbf G^{(i)}_k$ to eliminate $h^{(i)}_{k+1,k}$}
\State{$\mathbf H^{(i)} \gets \mathbf G^{(i)}_k \mathbf H^{(i)}$}
\Comment{apply $\mathbf G^{(i)}_k$ to $\mathbf H^{(i)}$}
\State{$\hat r^{(i)} \gets \mathbf G^{(i)}_k \hat r^{(i)}$}
\Comment{apply $\mathbf G^{(i)}_k$ to $\hat r^{(i)}$}
\EndParFor
\State{$k \gets k + 1$}
\EndWhile
\ParFor{$i = 1, \ldots , m$}
\State{Solve $\mathbf{R^{(k, i)}} y^{(i)} = \hat r^{(i)}$}
\Comment{$\mathbf{R^{(k, i)}}$: upper triangular part of $\mathbf H^{(i)}|_{k \times k}$ after Givens rotations}
\State{Update $x^{(i)} \gets x^{(i)} + \mathbf V^{(k, i)} y^{(i)}$}
\Comment{update $x^{(i)}$ with least-squares solution}
\EndParFor
\State{\Return $\mathbf x$}
\end{algorithmic}
\end{algorithm}

Vectors $q^{(k)} \in \mathbb K^n$ are now promoted to matrices $\mathbf q^{(k)} \in \mathbb K^{n \times m}$, and the $k$-th search direction of the $i$-th linear system is denoted by $q^{(k, i)} \in \mathbb K^n$.
The essential operations remain the matrix–matrix products in lines 2 and 8 of \cref{alg:pgmres_alg}.
The remaining computations can be parallelized via the \texttt{parfor} loops.
Each linear system constructs and operates exclusively on its own Krylov subspace.
Consequently, the right-hand side vectors $\mathbf b$, the solution vectors $\mathbf x$, and the search directions $\mathbf q^{(k)}$ naturally form dense matrices in $\mathbb K^{n \times m}$.
Other objects such as the Hessenberg matrices $\mathbf H^{(i)}$ and the transformed residuals $\hat r^{(i)}$ may be stored contiguously or not, although the former layout is preferable for efficient use of BLAS routines.

\subsubsection{Numerical examples}

\begin{figure}
\begin{subfigure}{.5\textwidth}
\includegraphics[width=\textwidth]{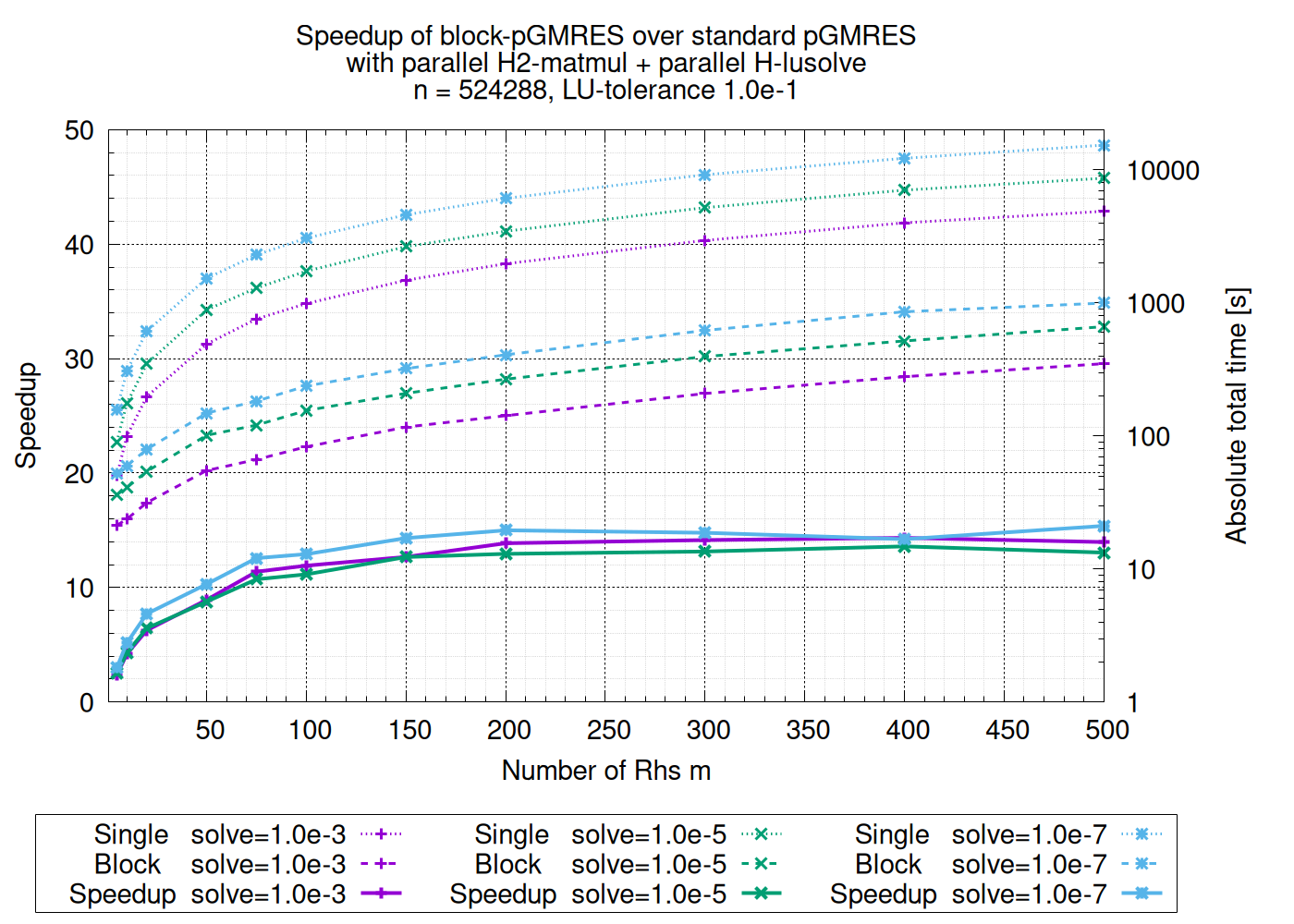}
\end{subfigure}%
\begin{subfigure}{.5\textwidth}
\includegraphics[width=\textwidth]{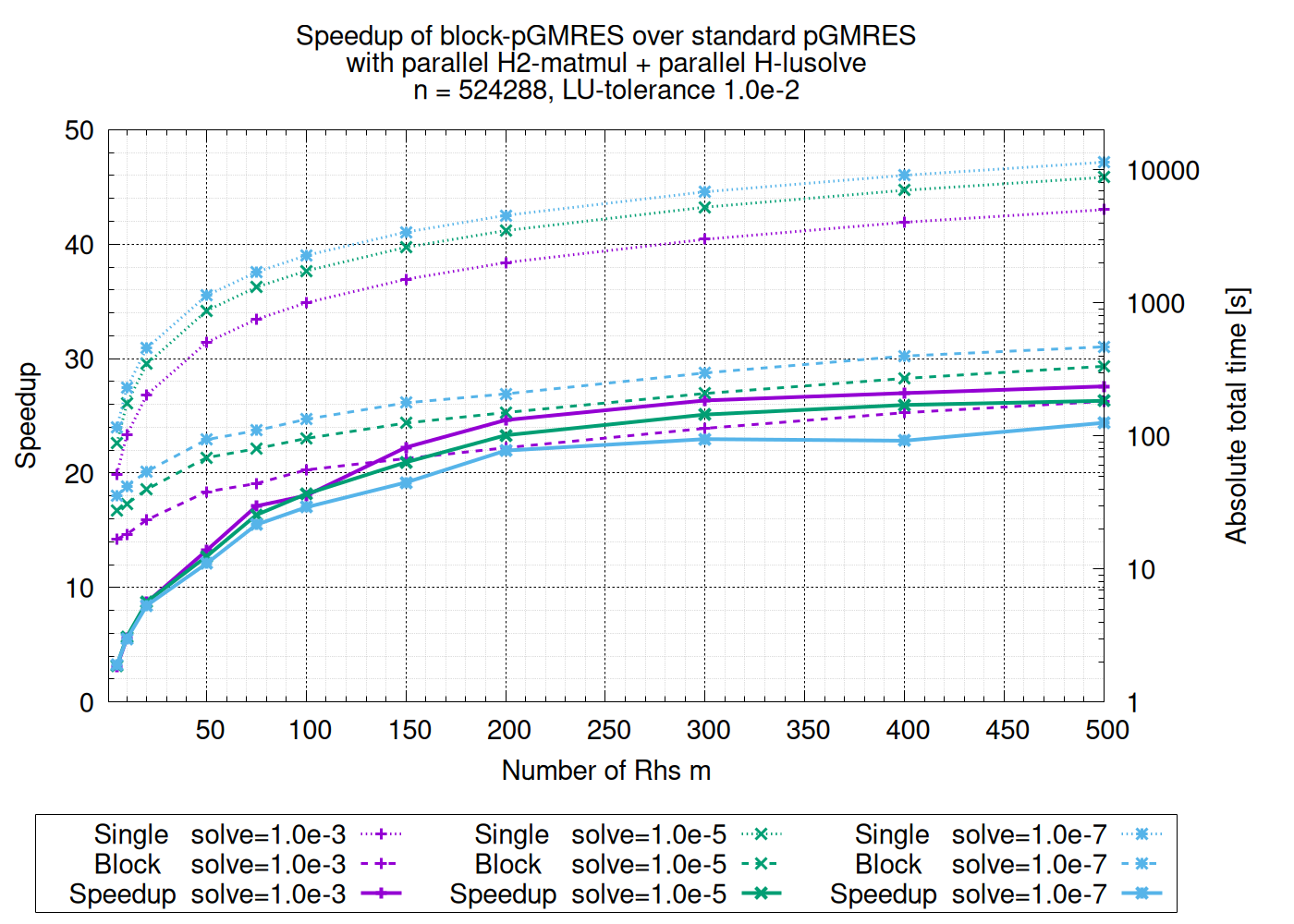}
\end{subfigure}%
\\
\begin{subfigure}{.5\textwidth}
\includegraphics[width=\textwidth]{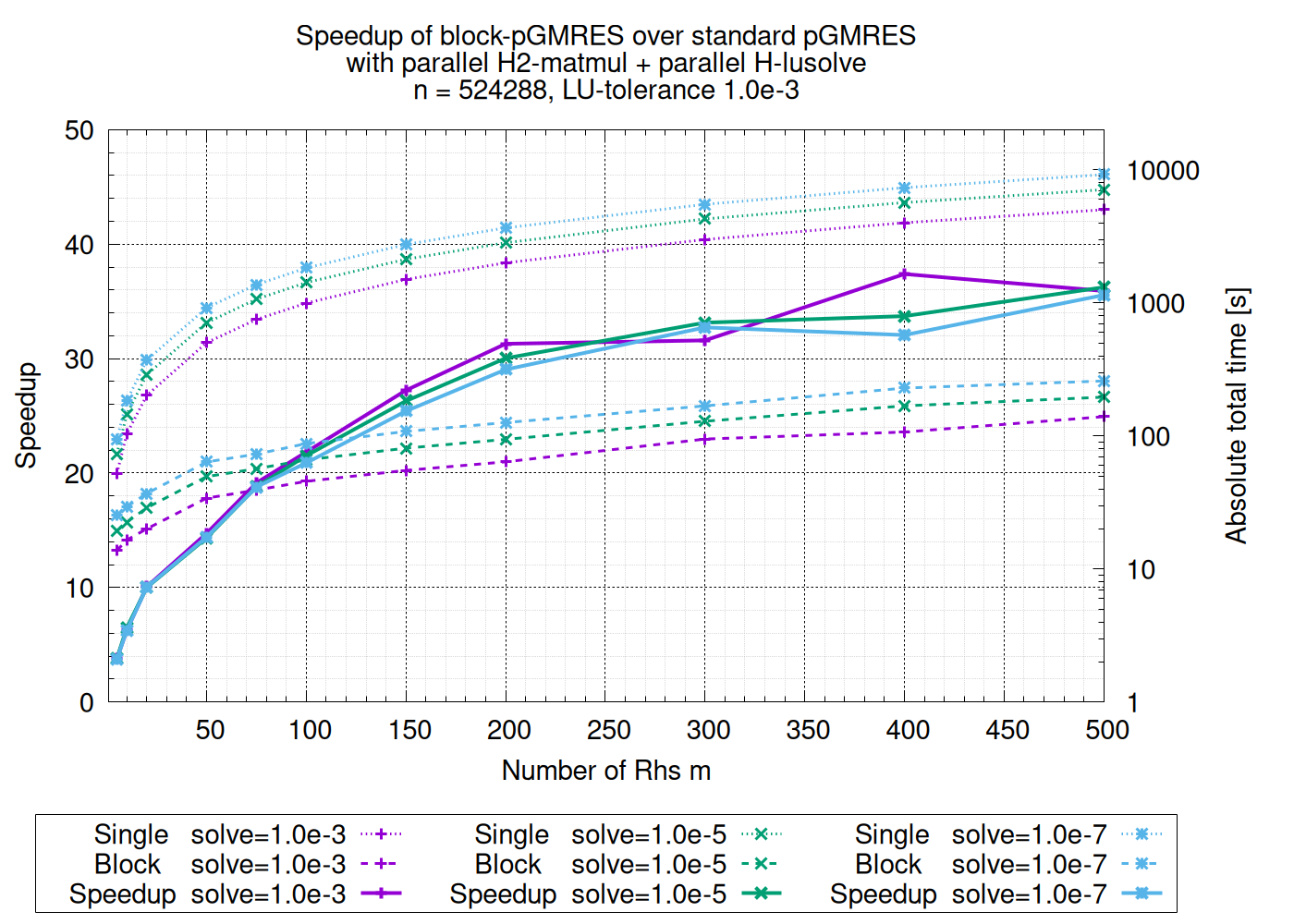}
\end{subfigure}%
\begin{subfigure}{.5\textwidth}
\includegraphics[width=\textwidth]{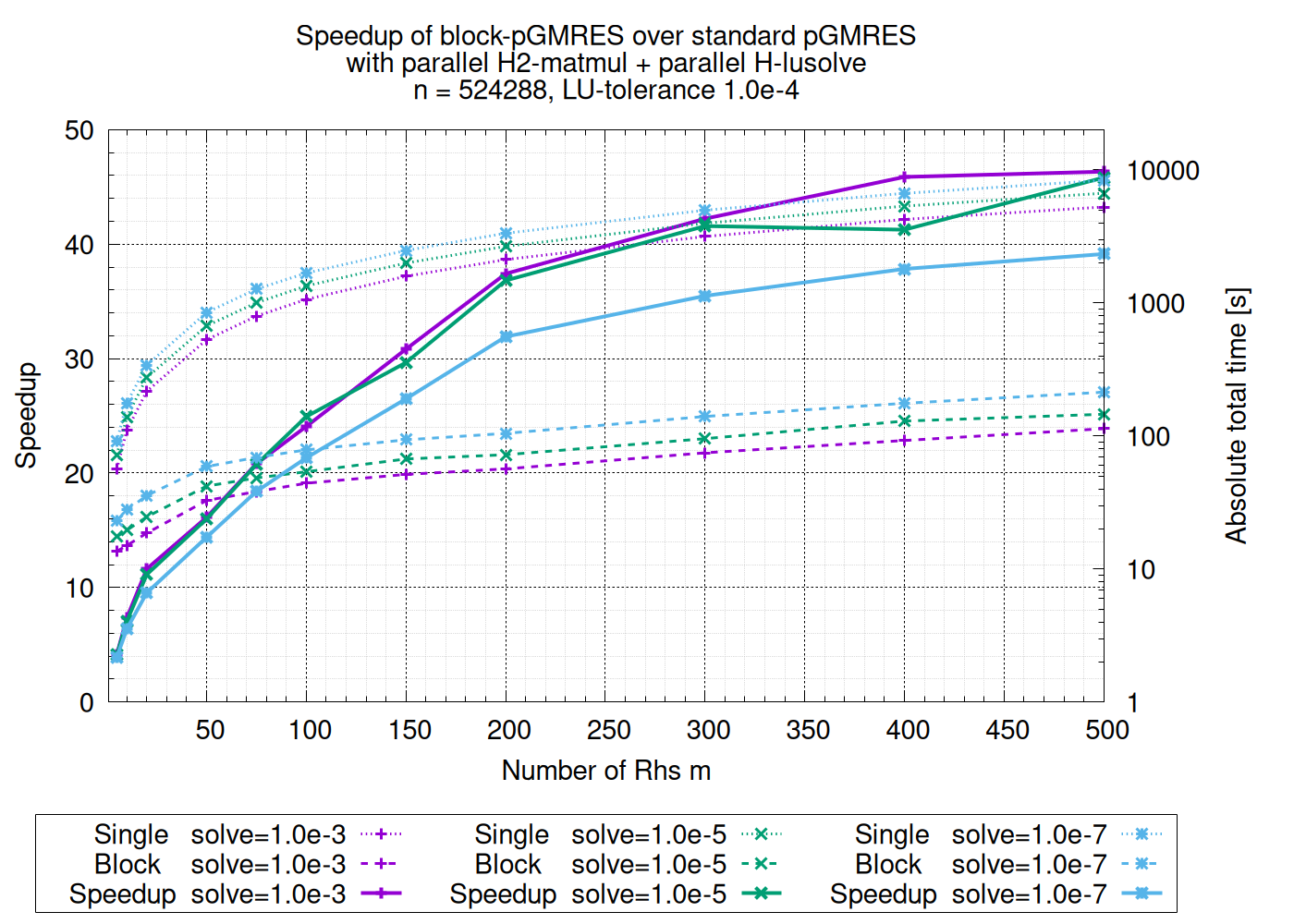}
\end{subfigure}%
\caption{Absolute runtimes for the standard pGMRES method and the block pGMRES method
with parallel $\mathcal H^2$-MVM and $\mathcal H^2$-MM, respectively, as well as the speedups achieved
by the latter compared to the former.
We vary the accuracy of the LU-preconditioner $\varepsilon_{\mathrm{dcp}}$ and the relative solver tolerance $\varepsilon_{\mathrm{slv}}$.}
\label{img:speedup_H2_pGMRES}
\end{figure}

As a motivating example for this method, we once again consider the single-layer operator
of the complex-valued Helmholtz equation with a low-frequency wavenumber of $\kappa = 1$ on the
unit sphere.
In this more general setup we choose the $\mathcal H$-LU-factorization as an efficient
and applicable preconditioner to the linear systems ~\cite{bebendorf2005hierarchical, carratala2019exploiting, boerm2022semiautomatic}.
As before, we present numerical results illustrating the performance of this approach
for different preconditioner accuracies and solver tolerances.
The corresponding absolute runtimes and speedups are shown in \cref{img:speedup_H2_pGMRES}.

The results appear very similar to those obtained with the CG method
and the highest speedup of 45 is also reached for $\varepsilon_{\mathrm{dcp}} = 10^{-4}$.
As the accuracy of the $\mathcal H$-LU preconditioner increases, the speedup of the method also improves. 
This indicates that more work is required for applying the preconditioner, which can
in turn be parallelized across the right-hand sides.
Since such parallelism is not available in the non-blocked version, the block-GMRES method
offers a clear efficiency advantage in this setting.

\begin{figure}
\begin{subfigure}{.5\textwidth}
\includegraphics[width=\textwidth]{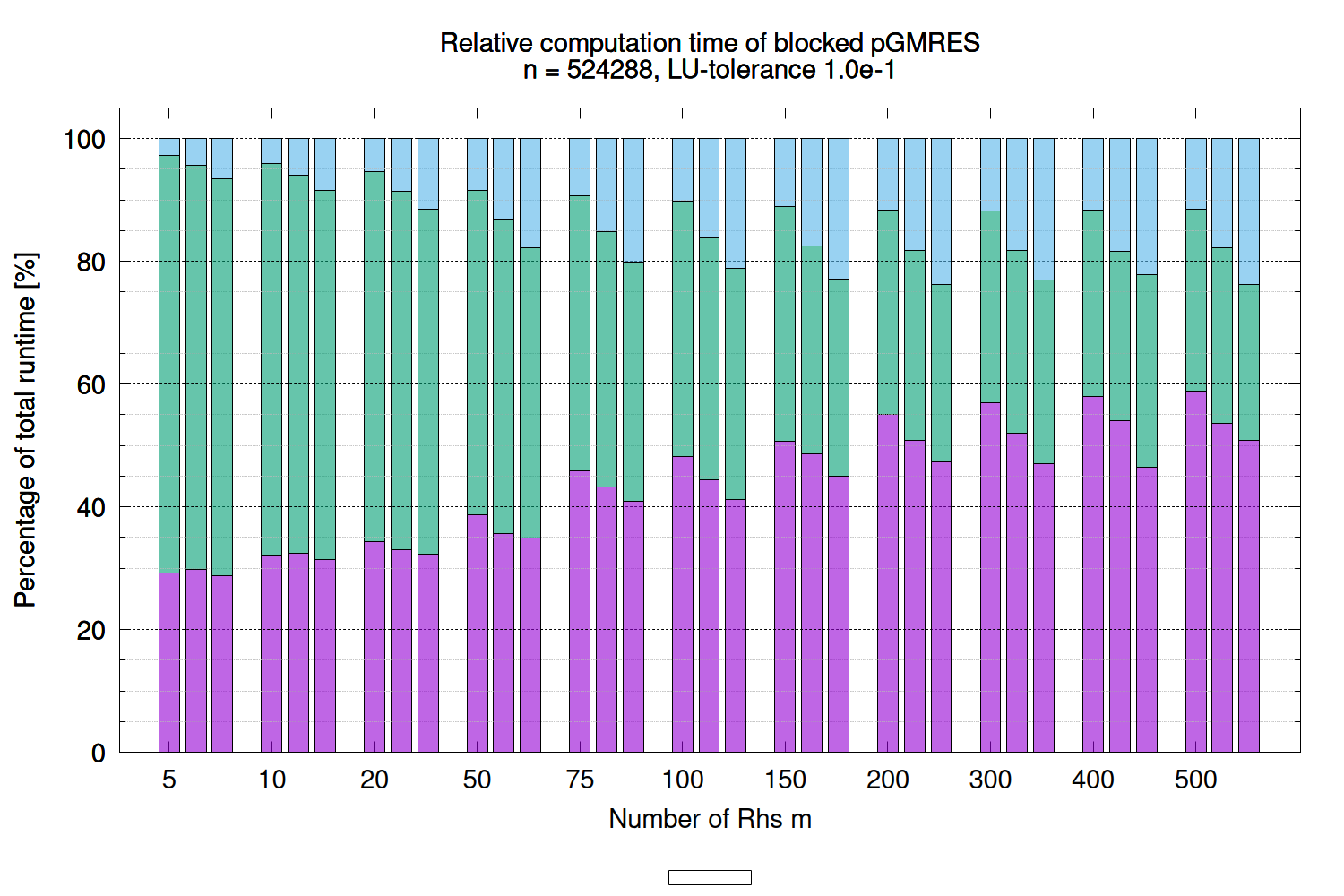}
\end{subfigure}%
\begin{subfigure}{.5\textwidth}
\includegraphics[width=\textwidth]{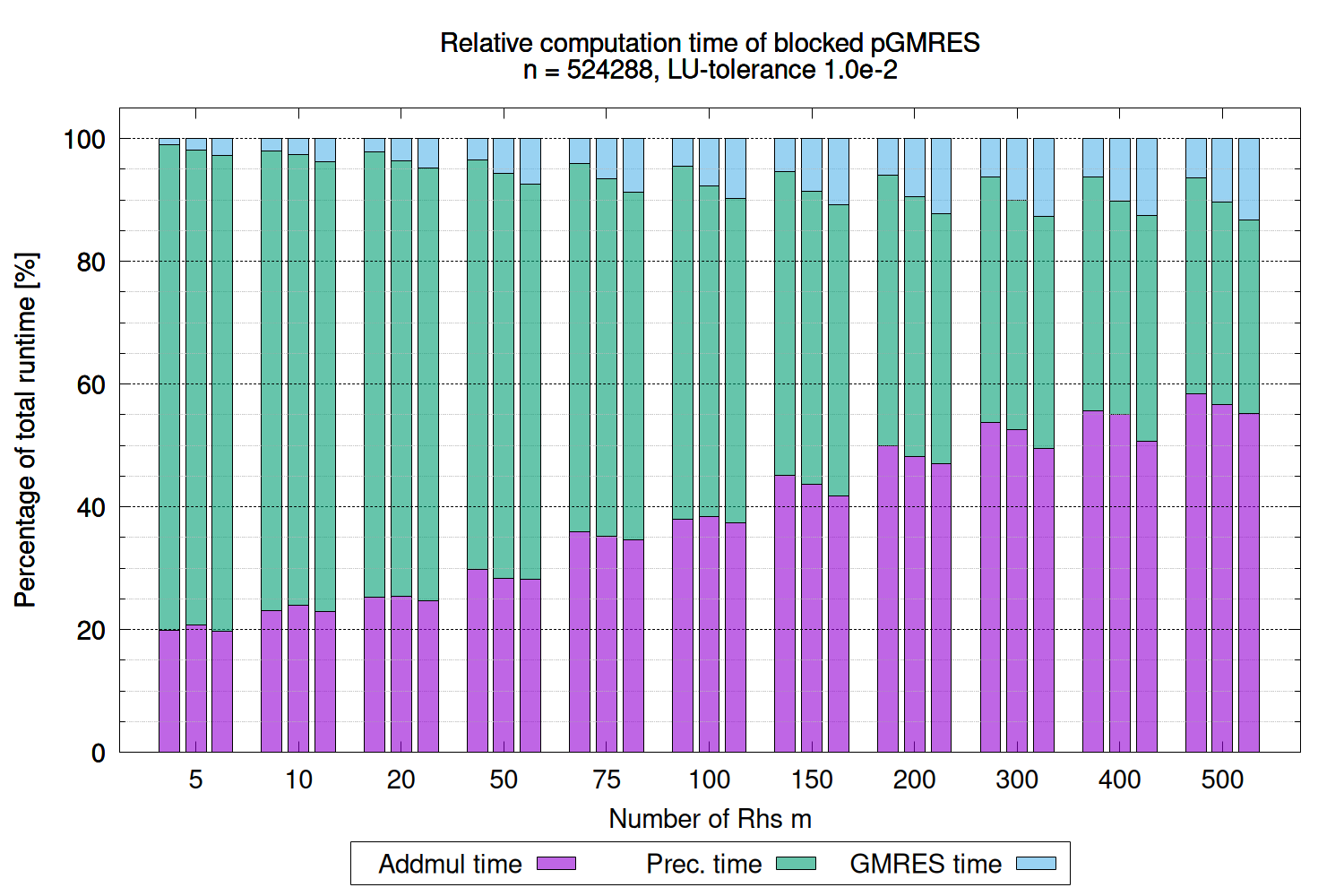}
\end{subfigure}%
\\
\begin{subfigure}{.5\textwidth}
\includegraphics[width=\textwidth]{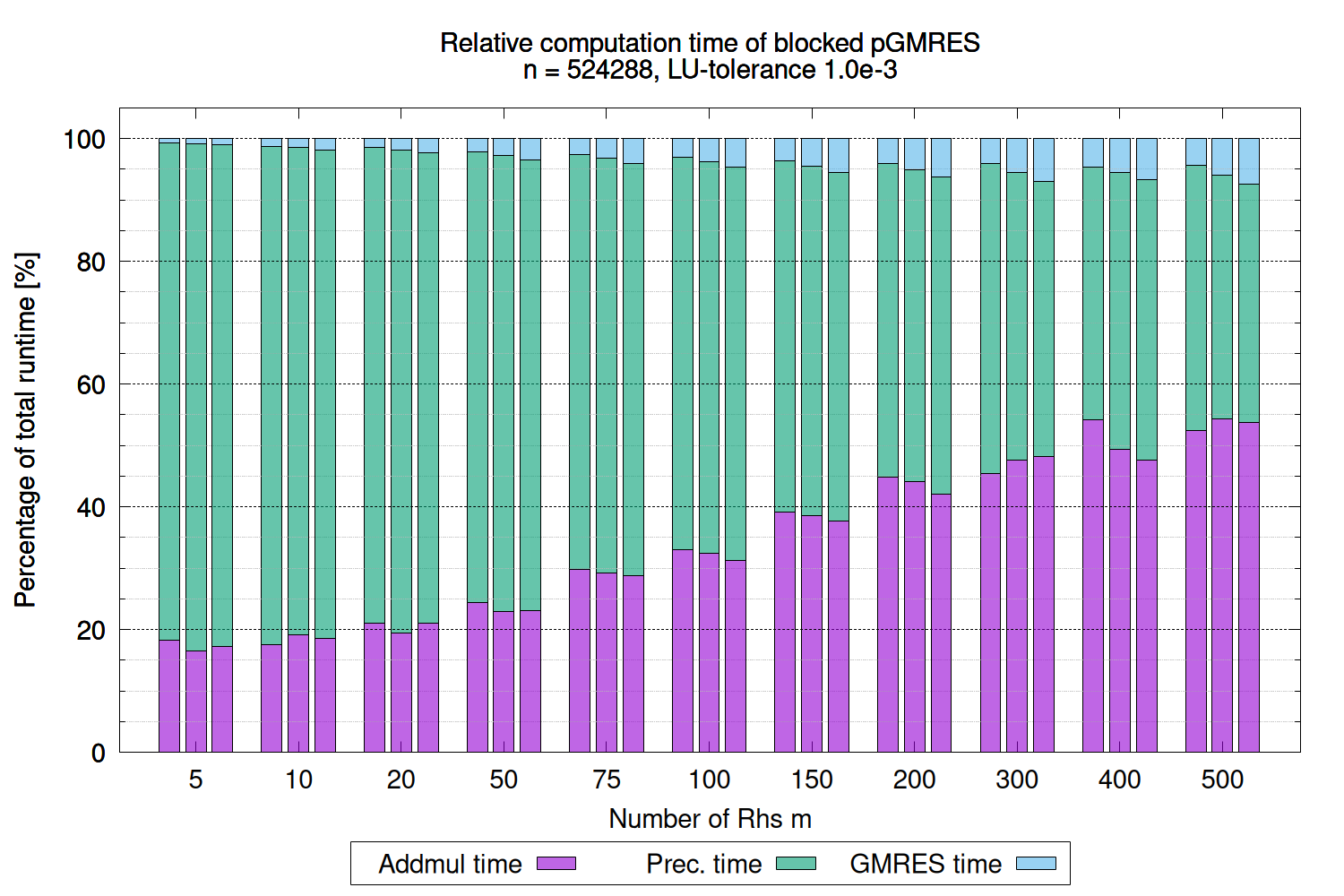}
\end{subfigure}%
\begin{subfigure}{.5\textwidth}
\includegraphics[width=\textwidth]{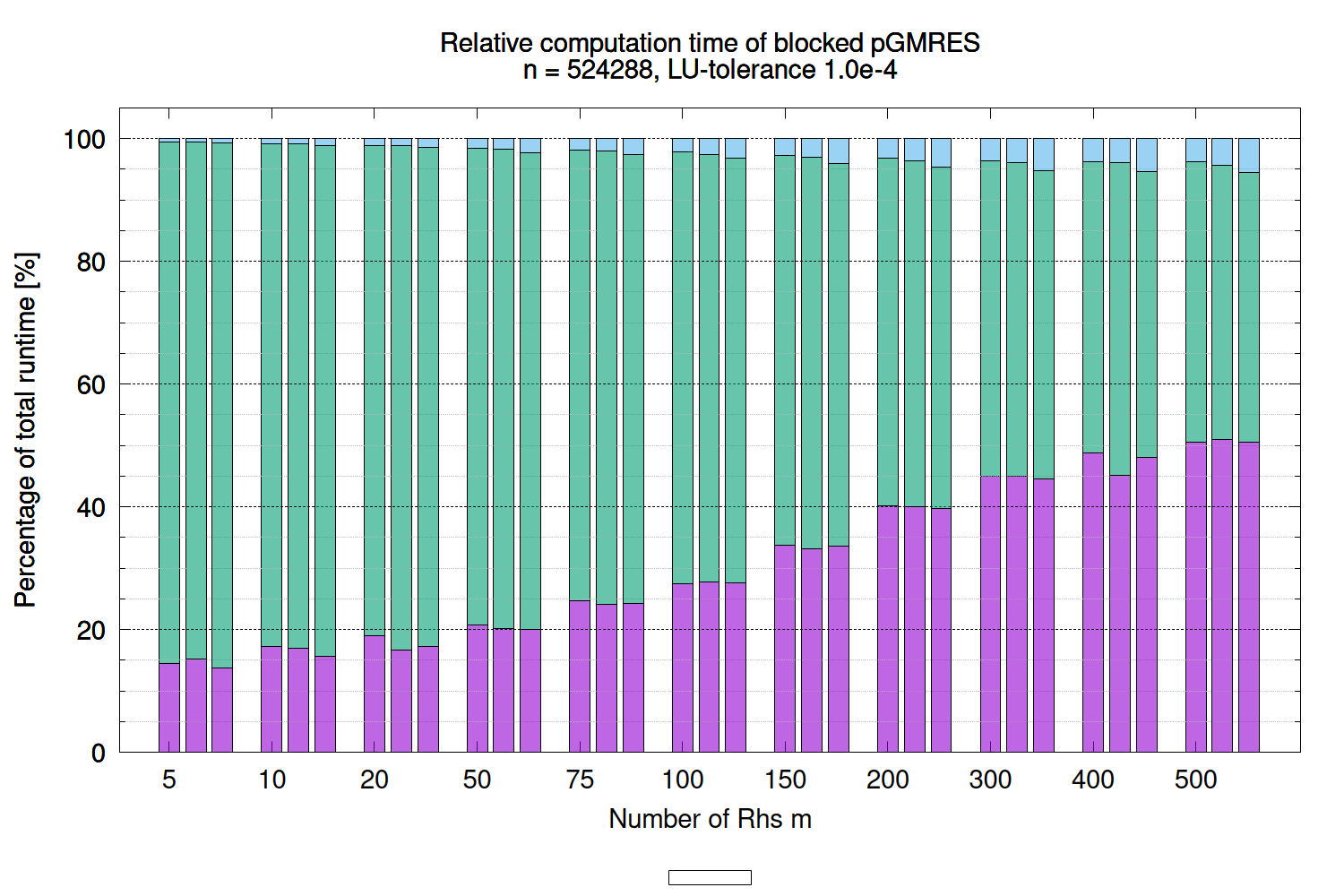}
\end{subfigure}%
\caption{Relative computing times of the solution process for \cref{alg:pgmres_alg}.
We vary the accuracy of the LU-preconditioner $\varepsilon_{\mathrm{dcp}}$ and the relative solver tolerance $\varepsilon_{\mathrm{slv}}$.
Within each group, the left bar corresponds to $\varepsilon_{\mathrm{slv}} = 10^{-3}$, the middle bar to $\varepsilon_{\mathrm{slv}} = 10^{-5}$, and the right bar to $\varepsilon_{\mathrm{slv}} = 10^{-7}$.}
\label{img:percentage_H2_pgmres}
\end{figure}

We now turn the readers attention to \cref{img:percentage_H2_pgmres}, which again
shows the relative runtime contributions of the matrix–matrix multiplication with the
system matrix $\mathbf A$, the application of the preconditioner $\mathbf M$, and the
remaining parts of the GMRES algorithm, denoted as ``GMRES time''.
The results are similar to those obtained for the CG method, but with some notable
differences.  
For small values of $m$, the preconditioner accounts for a larger fraction of the runtime
while the matrix–matrix multiplication is relatively inexpensive.  
As $m$ increases, this ratio shifts in favor of the matrix–matrix multiplication.  

A particularly noteworthy aspect is the proportion of time spent in the GMRES part itself,
which accounts for about $1$–$25\%$ of the total runtime.
Since each new iteration increases the cost of the Arnoldi process, one might expect this
percentage to rise for tighter solver tolerances, such as
$\varepsilon_{\mathrm{slv}} = 10^{-7}$, as shown in \cref{img:percentage_H2_pgmres}.
The better the preconditioner is constructed, the less time is required in the GMRES part,
since fewer iterations are needed.
In contrast, the time spent on evaluating the preconditioner $\mathbf M$
becomes increasingly expensive as $\varepsilon_{\mathrm{dcp}}$ decreases. 

\begin{remark}{Memory footprint}{}
\label{rem:memory_block_pgmres}
We now also consider the memory consumption of the GMRES and pGMRES methods.  
Here, the auxiliary vectors $r$, $\hat r$, $a$, and $q$ are required.  
For an efficient implementation of the Arnoldi process, one additionally needs to store the matrix 
$\mathbf H \in \mathbb K^{n \times \ell}$, assuming at most $\ell$ steps of the method are performed before a restart is triggered.  
Thus, in addition to the memory required for the $\mathcal H^2$-matrix itself, the dominant factor is the storage of $\mathbf H$, which requires $\mathcal O(n \ell)$
units of storage.  

The situation becomes more demanding in the block variant of the algorithm.  
In this case, a separate matrix $\mathbf H^{(i)}$ must be stored for each $i = 1, \ldots , m$, 
leading to a memory requirement of $\mathcal O(n \ell m)$ solely for the Hessenberg matrices.  
The auxiliary vectors add another $\mathcal O(n m)$ units of storage.  
Hence, the user must carefully monitor the total memory consumption.  
If $m$ becomes too large, it is advisable to process the right-hand sides in chunks of at most $m_{\mathrm{chunk}} < m$.  
\end{remark}

\newpage

\section{Conclusions}
\label{sec:conclusions}

We have presented an efficient strategy for parallelizing matrix–vector and 
matrix–matrix multiplication where the left factor is an $\mathcal H^2$-matrix. 
By eliminating race conditions in these operations, we showed that the 
matrix–vector product can fully utilize the available memory bandwidth and is thus 
optimal. Moreover, the systematic use of BLAS level-3 routines enables more 
cache-efficient operations, thereby exploiting the full potential of modern CPUs 
for matrix–matrix multiplication with $\mathcal H^2$-matrices. 
This operation is a key ingredient for solving multiple linear systems simultaneously 
within Krylov subspace methods. 
In combination with robust preconditioners such as 
$\mathcal H$-Cholesky or $\mathcal H$-LU factorizations, which can be constructed 
efficiently via $\mathcal H$-arithmetic, the approach proves highly effective. 
We demonstrated that combining parallel $\mathcal H^2$-matrix–matrix multiplication 
with iterative Krylov solvers outperforms solvers based solely on parallel 
matrix–vector multiplication by factors of 30–50 on a 24-core Intel processor.

Further performance improvements are possible. 
Parallel matrix–vector and matrix–matrix multiplication could be optimized in a 
NUMA-aware fashion or extended to distributed systems, for example using MPI. 
The scalability of preconditioner application, however, inherently depends on the 
availability of sufficiently many right-hand sides and eventually saturates with 
increasing core counts. 
One potential remedy is the use of task graphs to encode the operations required 
for evaluating the preconditioner. 
Prior work ~\cite{boerm2022semiautomatic, christophersen2025schnelle} has shown 
that such techniques yield notable speedups in the construction of hierarchical 
matrix factorizations and may also accelerate their evaluation. 
Nevertheless, these algorithms remain sequential at their core. 
An alternative would be to employ an explicit matrix inverse as preconditioner, 
which again reduces the problem to a parallelizable matrix–matrix product. 
However, this approach is well known to suffer from numerical instability and 
higher computational as well as memory costs.

Looking ahead, the presented techniques open up several promising directions. 
The efficient combination of hierarchical matrices with block Krylov methods is 
highly relevant for large-scale applications where many right-hand sides occur 
naturally, such as in frequency-domain simulations in computational electromagnetics, 
wave propagation, or uncertainty quantification. 
In these contexts, the ability to solve hundreds of systems simultaneously with 
high efficiency has the potential to significantly reduce simulation times and 
resource consumption.

From a methodological perspective, extending the proposed algorithms to hybrid 
CPU–GPU systems or fully distributed architectures represents an important next step. 
Here, the cache-aware BLAS level-3 formulations and the block-structured Krylov 
approach provide a strong foundation for leveraging the massive parallelism of 
modern hardware. 
Another promising direction lies in adaptive strategies that automatically balance 
the chunk size $m_{\mathrm{chunk}}$ against available memory and computational 
resources, making the methods even more flexible for practical use.

Ultimately, the synergy between hierarchical matrix arithmetic, robust 
preconditioning, and block Krylov solvers offers a pathway towards scalable and 
general-purpose iterative solvers. 
We believe that further exploration of these techniques will help bridge the gap 
between sophisticated numerical linear algebra and the practical demands of 
large-scale scientific and engineering applications.

\bibliographystyle{siamplain}
\bibliography{references}
\end{document}